\documentclass[12pt]{amsart}
\font\emailfont=cmtt10

\usepackage{amsmath,amsthm,amsfonts,amscd,flafter,epsf}
\usepackage[small,nohug,heads=littlevee]{diagrams} 
\diagramstyle[labelstyle=\scriptstyle] 

\newcommand{\rk}{\mathrm{rk}}

\newtheorem{theorem}{Theorem}[section]
\newtheorem{prop}[theorem]{Proposition}
\newtheorem{cor}[theorem]{Corollary}

\newtheorem{lemma}[theorem]{Lemma}

\newtheorem{defn}[theorem]{Definition}

\def\endproof{\relax\ifmmode\expandafter\endproofmath\else
  \unskip\nobreak\hfil\penalty50\hskip.75em\hbox{}\nobreak\hfil\bull
  {\parfillskip=0pt \finalhyphendemerits=0 \bigbreak}\fi}
\def\endproofmath$${\eqno\bull$$\bigbreak}
\def\bull{\vbox{\hrule\hbox{\vrule\kern3pt\vbox{\kern6pt}\kern3pt\vrule}\hrule}}

\newcommand{\R}{\mathbb{R}}

\newcommand{\Z}{\mathbb{Z}}

\newcommand{\OneHalf}{\frac{1}{2}}

\newcommand{\cm}{\cdot}

\newcommand{\ModSWfour}{\mathcal{M}}
\newcommand{\ModFlow}{\ModSWfour}

\newcommand\abuts\Rightarrow
\newcommand\Sym{\mathrm{Sym}}

\newcommand\relspinc{\underline{\spinc}}

\newcommand\x{\mathbf x}
\newcommand\w{\mathbf w}

\newcommand\y{\mathbf y}

\newcommand\ModSphere{\ModFlow\left({\mathbb S}\longrightarrow 
\Sym^{g-1}(\Sigma_{1})\times \Sym^2(\Sigma_{2})\right)}
\newcommand\ModSpheres\ModSphere
\newcommand\CF{CF}

\newcommand\CFa{\widehat{CF}}

\newcommand\CFm{\CF^-}

\newcommand\HFa{\widehat{HF}}

\newcommand\Mas{\mu}
\newcommand\UnparModSp{\widehat \ModSp}
\newcommand\UnparModFlow\UnparModSp
\newcommand\Mod\ModSp

\newcommand{\cald}{{\mathcal D}}

\newcommand{\spinc}{\mathfrak s}

\newcommand\ModMaps{\mathcal M}
\newcommand\ModSp\ModMaps

\newcommand\Ta{{\mathbb T}_{\alpha}}
\newcommand\Tb{{\mathbb T}_{\beta}}

\newcommand\alphas{\mbox{\boldmath$\alpha$}}

\newcommand\betas{\mbox{\boldmath$\beta$}}

\newcommand\Dual{\mathcal D}
\newcommand\Duality\Dual

\newcommand\Ko{\mathcal R}
\newcommand\Kp{\Knot_+}
\newcommand\Km{\Knot_-}

\newcommand\Kauff{\mathfrak K}

\newcommand\oLink{\vec{L}}
\newcommand\oKnot{\vec{K}}
\newcommand\dT{\underline{\partial}}
\newcommand\CFKmT{{\underline{CFK}}^-}
\newcommand\CFLmT{{\underline{CFL}}^-}

\newcommand\HFKmT{{\underline{\mathrm{HFK}}}^-}
\newcommand\CFKmTP{\CFKmT_{(2)}}
\newcommand\CFKmTI{\CFKmT_{(1)}}
\newcommand\CFKaT{\widehat{\underline{CFK}}}
\newcommand\Gen{\mathfrak{S}}

\newcommand\EmptyRect{\Rect^\circ}

\newcommand\CKm{CK^-}
\newcommand\CKa{\widehat{CK}}

\newcommand{\sign}{\mathcal S}

\newcommand\AlgM{\mathcal{A}}
\newcommand\RingM{\mathfrak R}

\newcommand\Rect{\mathrm{Rect}}
\newcommand\Zip{\mathcal Z}
\newcommand\Unzip{\mathcal U}
\newcommand\Sing{X}
\newcommand\As{\mathbb A}
\newcommand\Bs{\mathbb B}
\newcommand\Os{\mathbb O}
\newcommand\Ps{\mathbb P}
\newcommand\Xs{\mathbb X}

\newcommand\CnewM{{C}}

\newcommand\Knot{\mathcal K}
\newcommand\Link{\mathcal L}

\newcommand\Out{\mathrm{Out}}
\newcommand\In{\mathrm{In}}

\newcommand\orL{\vec{L}}

\newcommand\CFLa{\widehat{\mathrm CFL}}
\newcommand\CFLm{\mathrm{CFL}^-}
\newcommand\HFLm{\mathrm{HFL}^-}
\newcommand\HFLa{\widehat {\mathrm{HFL}}}

\newcommand\ws{\mathbf w}
\newcommand\zs{\mathbf z}

\newcommand\spincrel\relspinc

\newcommand\CFK{\mathrm{CFK}}
\newcommand\HFK{\mathrm{HFK}}

\newcommand\CFKa{\widehat\CFK}
\newcommand\CFKm{\CFK^-}

\newcommand\HFKa{\widehat\HFK}
\newcommand\HFKm{\HFK^-}

\newcommand\Smooth{\mathcal R}
\newcommand\Singularize{\mathcal X}

\title[{A cube of resolutions for knot Floer homology}]
{A cube of resolutions for knot Floer homology}

\author[Peter Ozsv{\'a}th]{Peter Ozsv\'ath}
\address{Department of
Mathematics, Columbia University, 
New York, NY 10027 \newline
\indent{\emailfont{petero@math.columbia.edu}}}
\thanks{PSO was supported by NSF grant number DMS-0505811}

\author[Zolt{\'a}n Szab{\'o}]{Zolt{\'a}n Szab{\'o}} 
\address{Department of
Mathematics, Princeton University, New Jersey 08544 \newline
\indent{\emailfont{szabo@math.princeton.edu}}}
\thanks{ZSz was supported by NSF grant number DMS-0406155}

\begin{document}
\begin{abstract}
  We develop a skein exact sequence for knot Floer homology, involving
  singular knots.  This leads to an explicit, algebraic description of
  knot Floer homology in terms of a braid projection of the knot.

\end{abstract}

\maketitle
\section{Introduction}
\label{sec:Introduction}

Knot Floer homology is an invariant for knots in $S^3$, defined using
Heegaard diagrams and holomorphic disks \cite{Knots},
\cite{RasmussenThesis}. One version gives a finitely generated,
bigraded Abelian group
$$\HFKa(\Knot)=\bigoplus_{d,s\in\Z}\HFKa_d(\Knot,s)$$
whose Euler
characteristic gives the symmetrized Alexander polynomial, in the
sense that
$$\sum_{d,s\in\Z} (-1)^d \rk ~\HFKa_d(\Knot,s) \cm
T^s=\Delta_\Knot(T).$$
The grading specified by $d$ is called the {\em
  Maslov grading}, and the one specified by $s$ is called the {\em
  Alexander grading}.  There is another variant of knot Floer
homology, $\HFKm(\Knot)$, which is also a bigraded module, admits an
action by the polynomial algebra $\Z[U]$ (where $U$ has bigrading
$(-2,-1)$).  Its bigraded Euler characteristic (over $\Z$) is given by
$\Delta_K(T)/(1-T)$.

In a recent paper~\cite{MOS}, knot Floer homology is calculated using
elementary combinatorial means, see also~\cite{MOST},
\cite{SarkarWang}. The aim of the present paper is to give a more
conceptual, algebraic description of the knot Floer homology groups of
a knot $\Knot$, given in terms of a braid presentation of $\Knot$.
This description leads to a purely algorithmic calculation of knot
Floer homology, which are both quite different in character, and also
logically independent, from the ones from~\cite{MOS}
and~\cite{SarkarWang}.

An important tool in this description is a skein exact sequence which
relates the knot Floer homology of a knot with its smoothing and its
singularization, established in Theorem~\ref{thm:SkeinExactSequence} below.

More precisely, recall the usual skein relation for the Alexander
polynomial, which states that the Alexander-Conway polynomial
$\Delta_\Knot(T)$ for oriented links is uniquely characterized by its
value on the unknot, together with the skein relation
\begin{equation}
  \label{eq:UsualSkein}
  \Delta_{\Knot_+}(T)-\Delta_{\Knot_-}(T)=(T^{1/2}-T^{-1/2})
  \cm \Delta_\Smooth(T),
\end{equation}
where here $\Knot_+$ and $\Knot_{-}$ are knots specified by two
projections which differ at a single crossing which is positive for
$\Knot_+$ and negative for $\Knot_-$; and $\Smooth$ denotes the
smoothing at this crossing.  Indeed, there is an exact triangle
generalizing the above skein relation
(cf.~\cite[Equation~\eqref{Knots:eq:SkeinExact}]{Knots}).  For that
exact triangle, of course, the notion of knot Floer homology is
extended to an invariant for oriented links.

A more useful exact triangle is derived in the present paper, which
involves singular links, compare also~\cite{KhovanovRozanskyII}. To this
end, recall that there is a skein relation involving a knot, its
smoothing at a given crossing $\Smooth$, and its singularization ${\mathcal X}$.  Specifically,
the Alexander polynomial for oriented links can be uniquely extended
to singular oriented links (i.e. with double-point singularities) to satisfy the
skein relations
\begin{eqnarray}
 \Delta_{\Knot_+}= \Delta_{\mathcal X}+T^{\OneHalf}\Delta_\Smooth
&{\text{and}}&
\Delta_{\Knot_-}= \Delta_{\mathcal X}+T^{-\OneHalf}\Delta_\Smooth
 \label{eq:SingularSkein}
\end{eqnarray}
Of course, Equation~\eqref{eq:UsualSkein} follows quickly from the
above two equations.

For the exact triangle corresponding to
Equation~\eqref{eq:SingularSkein}, we must appropriately interpret the 
invariant for links and for singular knots.

Knot Floer homology is extended to an invariant $\HFLm(\orL)$ for
oriented links in~\cite{Links}. The invariant there is a module over
the polynomial algebra $\Z[U_1,...,U_\ell]$, where the $U_i$ are
variables which are in one-to-one correspondence with the link's
components. This invariant comes equipped with a Maslov grading, and
$\ell$ different Alexander gradings. For our present purposes, we need
only consider a collapsed Alexander grading, gotten by adding up all
$\ell$ of the Alexander gradings, so that we continue to think of the
link Floer homology as bigraded.

Furthermore, knot Floer homology is extended to (oriented) singular
links in~\cite{SingLink}. For a knot with a single singular point this
is the homology group of a chain complex $\CKa({\mathcal X})$ over the
ring $\Z[U_a,U_b,U_c,U_d]$, where the indeterminates belong to the
four edges at the singular point ($a$ and $b$ point out and $c$ and
$d$ point in).

In this paper, we establish the following exact sequence generalizing
Equation~\eqref{eq:SingularSkein}:

\begin{theorem}
  \label{thm:SkeinExactSequenceIntro}
  Fix a projection of a knot with a distinguished crossing $p$.
  Let $\Knot_+$ denote the knot with a positive crossing at $p$,
  $\Knot_-$ denote the knot with a negative crossing at $p$,
  $\Singularize$ its singularization at $p$,
  and $\Smooth$ its smoothing at $p$. 
  Then, there are long exact sequences
  $$
  \begin{CD}
    ...
    @>>>\HFKm(\Knot_+) 
    @>{\phi_+}>>
    H_*(\frac{\CFKm(\Singularize)}{(U_a+U_b- U_c- U_d)}) 
    @>{\psi_+}>>
    T^{\OneHalf} \cm \HFLm(\Smooth)@>{\delta_+}>>...,
  \end{CD}
  $$
  and 
  $$
  \begin{CD}
    ...
    @>>>
    T^{-\OneHalf}\cm \HFLm(\Smooth)
    @>{\phi_-}>> 
    H_*(\frac{\CFKm(\Singularize)}{(U_a+U_b-U_c-U_d)})
    @>{\psi_-}>>
    \HFKm(\Knot_-)
    @>{\delta_-}>> ...
  \end{CD}
  $$
  where the maps $\phi_+$, $\delta_+$, $\psi_-$, and $\delta_-$
  preserve Maslov gradings, whereas the maps $\phi_+$ and $\phi_-$
  drop it by one. Moreover, all maps preserve Alexander gradings, with
  the understanding that $T^{\pm \OneHalf}\cm \HFLm(\Smooth)$ denotes
  $\HFLm({\mathcal R})$, thought of as bigraded, and whose (combined)
  Alexander grading is shifted up by $\pm \OneHalf$.
\end{theorem}

A related version holds also for $\HFKa$, cf.
Theorem~\ref{thm:SkeinExactSequence} below.

Iterating the above exact sequence (on the chain level, cf.
Section~\ref{sec:ExactSeq}), we arrive at a description of the knot
Floer homology groups of an arbitrary knot in terms of the knot Floer
homology groups of fully singular knots. Combining this with an explicit
calculation of knot Floer homology for fully singular knots, we obtain a 
concrete description of knot Floer homology in general.  Before giving
the description, we explain some of the objects which go into this
description.

Fix a projection $\Knot$ of an oriented knot $K$, which is in braid
position, let $c(\Knot)$ denote its crossings, and let $n$ denote the
number of crossings. We distinguish one of the leftmost edges.
We call such a diagram a {\em decorated braid position} for a knot.
Let $E=\{e_0,...,e_{2n}\}$ denote the set of edges in the diagram, 
where we view
the distinguished edge is subdivided in two, so that $e_0$ is the
second segment in the circular ordering induced by the orientation of
$K$. 

Each crossing of $\Knot$ can be either singularized or smoothed. A
{\em singularization} of $\Knot$ at a crossing $p$ is a singular link
with a double-point at $p$.  A {\em smoothing} of $\Knot$ at $p$ is a
link which is given the oriented resolution of $\Knot$ at $p$, i.e.
the crossing is replaced by a diagram with one fewer crossing.
Assigning smoothings or singularizations at each crossing of $\Knot$,
we obtain a planar graph $G$ where each vertex has valence either four
(if it corresponds to a singularized crossing) or two (if it
corresponds to a smoothed crossing; or more precisely, each smoothed
crossing gives rise to two valence two vertices). If $S$ is a singular
link obtained from $\Knot$ by singularizing or smoothing each crossing
of $\Knot$, we call $S$ an associated {\em complete resolution} of
$\Knot$. 

Let $\RingM$ denote the algebra $\Z[t,U_0,...,U_{2n}]$.  The variables
$U_i$ correspond to the edges of our singular link, and we introduce
another indeterminate $t$ as well. We think of this ring as graded by
an Alexander grading, where $U_i$ has Alexander grading $-1$, and
$\Z[t]$ is supported in Alexander grading equal to zero. Sometimes,
we will find it convenient also to complete $\RingM$, to obtain a new
ring $\RingM[[t]]$ of formal power series in $t$, with coefficients in 
the polynomial algebra $\Z[U_0,...,U_{2n}]$.

\begin{figure}[ht]
\mbox{\vbox{\epsfbox{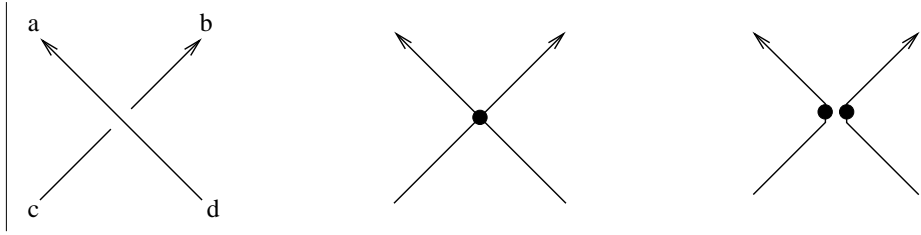}}}
\caption{\label{fig:ResNegCross}
{\bf{Resolutions of crossings.}}
        We show the singularization (middle) and the smoothing (right)
        of the negative crossing (left).}
\end{figure}

We define now an $\RingM$-algebra $\AlgM(S)$ associated to some
fully singularized link $S$ for $\Knot$. This algebra
$\AlgM(S)$ is a quotient of
$\RingM$ by an ideal which we describe presently.  At each crossing $p$ of
$\Knot$, label the edges $\{a,b,c,d\}$, where $a$ and $b$ are out-going,
$c$ and $d$ are in-coming, and $a$ resp. $c$ is to the left of $b$
resp. $d$, cf. Figure~\ref{fig:ResNegCross}.

If $p$ is a singularized point in $S$, we introduce the relations
\begin{eqnarray}
t\cm U_a + t\cm U_b &=& U_c + U_d  \label{eq:LinRel} \\
t^2\cm U_a\cm  U_b &=& U_c\cm  U_d. \label{eq:QuadRel}
\end{eqnarray}
If $p$ is a smoothed point in $S$, we introduce the relations
\begin{eqnarray}
t\cm U_a &=& U_c \label{eq:ResRel1}\\
t\cm U_b &=& U_d. \label{eq:ResRel2}
\end{eqnarray}
Sometimes, we label these variables $U_a^{(p)}$, $U_b^{(p)}$, $U_c^{(p)}$,
$U_d^{(p)}$, when it is necessary to recall the singular point $p$.

We introduce further relations associated to each collection $W$ of
vertices in the graph associated to $S$.  The {\em weight} of $W$,
denoted $|W|$, is twice the number of singular vertices plus the
number of smoothed vertices in $W$. Let $\Out(W)$ denote the outgoing
edges of $W$, let $\In(W)$ denote the incoming ones, and let $W^c$
denote the complement of $W$.  For each collection $W$ of vertices
which does not contain the distinguished vertex, we introduce the
following relation
\begin{equation}
  \label{eq:GenRel}
  t^{|W|}\cm \prod_{e\in \Out(W)\cap \In(W^c)} U_e
 = \prod_{e\in \In(W)\cap \Out(W^c)} U_e
\end{equation}
(These relations generalize
Equations~\eqref{eq:QuadRel}, \eqref{eq:ResRel1}, and \eqref{eq:ResRel2}
in a straightforward manner: in each case, we consider the case where $W$
consists of a single vertex.)

The algebra $\AlgM(S)$ is obtained by dividing out $\RingM$ by the
relations in Equations~\eqref{eq:LinRel}, \eqref{eq:QuadRel},
\eqref{eq:ResRel1}, \eqref{eq:ResRel2}, and \eqref{eq:GenRel} above.
Note that since the relations are all homogeneous, the algebra
$\AlgM(S)$ inherits an Alexander grading from $\RingM$.  We define the
symmetrized Alexander grading on $\AlgM(S)$ by the formula
$$A=A_0 + \frac{\#(\sigma(S)) - b + 1}{2},$$
where $\sigma(S)$ denotes the
number of double points of $S$,
where $A_0$ is the original Alexander grading and $b$
denotes the braid index of $S$. 
We also define the {\em Maslov grading} on $\AlgM(S)$ to be twice the
symmetrized Alexander grading.

Note that when $S$ is disconnected,
the completion $\AlgM(S)[[t]]$ is trivial.  Specifically, taking a set
non-empty $W$ of vertices in a component of $S$ which does not contain
the distinguished edge, we see that both $\In(W)$ and $\Out(W)$ are
empty, and hence we obtain the relation $t^{|W|}-1=0$.

Let $\Smooth$ and $\Singularize$ be two singular links whose
projections which differ at only one crossing $p$, where ${\mathcal
  R}$ is smoothed at $p$, and $\Singularize$ is singularized at $p$.
We define $\RingM$-module homomorphisms
\begin{eqnarray*}
u_p\colon \AlgM(\Singularize)\longrightarrow \AlgM(\Smooth)
&{\text{and}}&
z_p\colon \AlgM(\Smooth)\longrightarrow \AlgM(\Singularize),
\end{eqnarray*}
the {\em unzip} and {\em zip} homomorphisms as follows.

The unzip homomorphism is the natural ring map, obtained by observing
that $\AlgM(\Smooth)$ is a quotient of $\AlgM(\Singularize)$. The
zip homomorphism is induced from multiplication by $t\cm U_a-U_d$. It
is an exercise in the relation from 
Equation~\eqref{eq:GenRel} to see that in fact these maps give
well-defined $\RingM$-module maps.

Given a decorated knot projection $\Knot$, we organize its set of
associated singular links in a ``cube of resolutions'',
compare~\cite{KhovanovRozansky}, \cite{KhovanovRozanskyII}. If $p$ is
a positive crossing, its $0$-resolution is singular at $p$, and its
$1$-resolution is its smoothing at $p$. If $p$ is a negative crossing,
then its $0$-resolution is its smoothing at $p$, while its
$1$-resolution its singularization at $p$.  With the above conventions, given
an assignment $I\colon c(\Knot)\longrightarrow
\{0,1\}$, we can form the associated singular link $X_{I}(\Knot)$.

We define an Abelian group 
$$\CnewM(\Knot)=\bigoplus_{I\colon c(\Knot)\longrightarrow \{0,1\}}
\AlgM(X_I(\Knot)).$$
The group $\CnewM(\Knot)$ is bigraded, i.e.
$$\CnewM(\Knot)\cong \bigoplus_{m,s\in\Z}\CnewM_{m,s}(\Knot).$$
The
two gradings are called the {\em Maslov grading} and the {\em
  Alexander grading}.  The {\em Maslov grading} on
$\AlgM(X_I(K))\subset \CnewM(\Knot)$ is induced from the
internal  Maslov grading on $\AlgM(X_I(K))$.
The {\em Alexander  grading} on
$\AlgM(X_I(K))\subset \CnewM(\Knot)$ is defined to be 
\begin{equation}
  \label{eq:RenormalizedAlexander}
  A'=A+ \OneHalf \left(-N+\sum_{p\in c(\Knot)} I(p)\right),
\end{equation}
where $A$ is the internal Alexander grading on $\AlgM(X_I(\Knot))$,
and $N$ denotes the number of negative
crossings in $L$.

We think of the index set of resolutions $I\colon
c(\Knot)\longrightarrow \{0,1\}$ as points in a hypercube (given the
obvious partial ordering) If two
vertices $I$ and $J$ of this hypercube are connected by an edge, then
$X_{I}(\Knot)$ is obtained from $X_{J}(\Knot)$ by a zip or unzip move.

Fix a pair of resolutions $I, J \colon c(\Knot)\longrightarrow \{0,1\}$
where $I$ and $J$ differ in a single place, and $I<J$, i.e.
there is some $p\in c(\Knot)$ with $I(p)=0$, $J(p)=1$, and $I(q)=J(q)$ 
for all $q\neq p$. Then, we define a map
$$D_{I<J}\colon \AlgM(X_I(\Knot)) \longrightarrow 
\AlgM(X_J(\Knot)),$$
as follows. Note that $X_{J}(\Knot)$ is obtained from $X_{I}(\Knot)$ by a
single zip or unzip move, and the map $D_{I<J}$ is
given by the zip or unzip homomorphism defined above. 
We can add these maps up to get a map
$$D\colon \CnewM(\Knot) \longrightarrow \CnewM(\Knot),$$
which is the
sum of $D_{I<J}$ over all $I,J\colon c(\Knot)\longrightarrow \{0,1\}$
as above.

Suppose that $I<K$ differ in two places. Then there are two 
distinct intermediate resolutions $J$ and $J'$ with $I<J<K$ and
$I<J'<K$. Note that
$$D_{J<K}\circ D_{I<J}
=D_{J'<K}\circ D_{I<J'}.$$
Thus, if we fix a map $\epsilon$ from edges $I<J$ in the hypercube to $\{\pm 1\}$
with the property that for all $(I,J,K)$ as above,
$\epsilon(J<K)\epsilon(I<J)=-\epsilon(J'<K)\epsilon(I<J')$,
we have that
the endomorphism
$$D=\sum_{I<J}\epsilon(I<J) D_{I<J}$$ is
satisfies $D^2=0$.

Observe also that $D$ drops Maslov grading by one and preserves Alexander
grading.

Let $\Z[t^{-1},t]]$ denote the ring of formal power series in $t$,
with integral coefficients, and with finitely many non-zero terms
consisting of $t$ raised to negative powers.  Similarly, let
$$\CnewM(K)[t^{-1},t]]=\bigoplus_{m,s\in\Z} \CnewM_{m,s}(K)[[t]]$$
denote the complex constructed from $\CnewM_{m,s}(K)[[t]]$, the
$\Z[t^{-1},t]]$-module consisting of formal power series
$$\sum_{n\in\Z} c_n t^n,$$
with $c_n\in \CnewM_{m,s}(K)$, and where $c_n=0$
for all sufficiently small $n$. The differential on $\CnewM(K)$ naturally
induces one on $\CnewM(K)[t^{-1},t]]$.

\begin{theorem}
  \label{thm:Calculate} Fix a projection of a knot $\Knot$ in $S^3$ as
  above.  Then, we have identifications \begin{eqnarray*}
    \HFKa(K)\otimes \Z[t^{-1},t]] &\cong&
    H_*(\CnewM(\Knot)[t^{-1},t]]/U_0=0) \\
    \HFKm(K)\otimes_{\Z} \Z[t^{-1},t]] &\cong& 
    H_*(\CnewM(\Knot)[t^{-1},t]]). \\
  \end{eqnarray*}
\end{theorem}

As mentioned earlier, Theorem~\ref{thm:Calculate} proved using the
following two ingredients.  First, we calculate the homology groups of
completely singular knots, identifying these groups with the algebra
described above.  Second, we establish a suitably general
version of Theorem~\ref{thm:SkeinExactSequenceIntro} relating the knot
Floer homology groups of a knot, its smoothing, and its
singularization.  Iterating the skein exact sequence at each crossing,
we obtain the ``cube of resolutions'' calculating $\HFKm$ above.

Clearly, specializing the above constructions to $t=1$, one obtains a
picture very closely connected with the HOMFLY homology of Khovanov
and Rozansky~\cite{KhovanovRozanskyII}. We hope to return to this point in a future paper.

This paper is organized as follows.

In Section~\ref{sec:KnotFloerHomology}, we recall briefly the
construction of knot Floer homology, and also the version for singular
knots~\cite{SingLink}.  In Section~\ref{sec:CalcSingLink}, we
calculate the knot Floer homology for purely singular knots,
identifying them with the algebras defined above. In
Section~\ref{sec:ExactSeq}, we
establish Theorem~\ref{thm:SkeinExactSequenceIntro}, and also a
generalization which gives a cube of resolutions description of knot
Floer homology.  This immediately gives a spectral sequence whose
$E_1$ term consists of Floer homology groups for the various
singularizations of the original knot.  In Section~\ref{sec:Proof},
the $d_1$ differential of the above spectral sequence is calculated,
and it is shown that the spectral sequence collapses after this stage,
giving the proof of Theorem~\ref{thm:Calculate}.
Finally, in Section~\ref{sec:Grids}, we
explain the singular exact triangle from the point of view of grid
diagrams.

\subsection{Acknowledgements}
The authors wish to thank Bojan Gornik, Mikhail Khovanov, Robert
Lipshitz, Ciprian Manolescu, Andr{\'a}s Stipsicz, and Dylan Thurston
for interesting discussions. Moreover, this paper grew out of earlier
work in collaboration with Andr{\'a}s, see~\cite{SingLink}.

\section{Knot Floer homology}
\label{sec:KnotFloerHomology}

Knot Floer homology is a bigraded Abelian group associated to a knot
in $S^3$, cf.~\cite{Knots}, \cite{RasmussenThesis}.  We will briefly
sketch this construction, and refer the reader to the above sources
for more details.

Let $\Sigma$ be a surface of genus $g$, let
$\alphas=\{\alpha_1,...,\alpha_{g+n-1}\}$ be a collection of pairwise
disjoint, embedded closed curves in $\Sigma$ which span a
$g$-dimensional subspace of $H_1(\Sigma)$. This specifies a handlebody
$U_\alpha$ with boundary $\Sigma$. Moreover,
$\alpha_1\cup...\cup\alpha_{g+n-1}$ divides $\Sigma$ into $n$
components, which we label
$$\Sigma-\alpha_1-...-\alpha_{g+n-1}={\mathfrak A}_1\coprod...\coprod {\mathfrak A}_n.$$
Fix another such collection of
curves $\betas=\{\beta_1,...,\beta_{g+n-1}\}$, giving another
handlebody $U_\beta$. 
Write
$$\Sigma-\beta_1-...-\beta_{g+n-1}={\mathfrak B}_1\coprod...\coprod
{\mathfrak B}_n.$$
Let $Y$ be the three-manifold specified by the
Heegaard decomposition specified by the handlebodies $U_\alpha$ and
$U_\beta$.  Choose collections of disjoint points
$\Os=\{O_1,...,O_n\}$ and $\Xs=\{X_1,...,X_n\}$, which are distributed
so that each region ${\mathfrak A}_i$ and ${\mathfrak B}_i$ contains
exactly one of the points in $\Os$ and also exactly one of the points
in $\Xs$.  We can use the points $\Os$ and $\Xs$ to construct an
oriented, embedded one-manifold $\oLink$ in $Y$ by the following
procedure. Let $\xi_i$ be an arc connecting the point in $\Xs\cap
{\mathfrak A}_i$ with the point in $\Os\cap {\mathfrak A}_i$, and let
$\xi_i'$ be its pushoff into $U_\alpha$ i.e. the endpoints of $\xi_i'$
coincide with those of $\xi_i$ (and lie on $\Sigma$), whereas its
interior is an arc in the interior of $U_\alpha$. The arc is endowed
with an orientation, as a path from an element of $\Xs$ to an element
of $\Os$.  Similarly, let $\eta_i$ be an arc connecting $\Os\cap
{\mathfrak B}_i$ to $\Xs\cap {\mathfrak B}_i$, and $\eta_i'$ be its
pushoff into $U_\beta$. Putting together the
$\xi_i'$ and $\eta_i'$, we obtain an oriented link $\oLink$ in $Y$.

\begin{defn}
The data $(\Sigma,\alphas,\betas,\Os,\Xs)$ is called a 
{\em pointed Heegaard diagram compatible with the oriented link
$\oLink\subset Y$}.
\end{defn}

An oriented link in a closed three-manifold $Y$ always admits a
compatible pointed Heegaard diagram. In this article, we will restrict
attention to the case where the ambient three-manifold $Y$ is $S^3$.

Consider the $g+n-1$-fold symmetric product of the surface $\Sigma$,
$\Sym^{g+n-1}(\Sigma)$.  This is equipped with a pair of tori
\begin{eqnarray*}
\Ta=\alpha_1\times...\times\alpha_{g+n-1} 
&{\text{and}}&
\Tb=\beta_1\times...\times\beta_{g+n-1}.
\end{eqnarray*}
Let $\Gen$ denote the set of intersection points $\Ta\cap\Tb\subset\Sym^{g+n-1}(\Sigma)$.

Let $\CFKm(\oLink)$ be the free module over $\Z[U_1,...,U_n]$
generated by elements of $\Gen$, where here the $\{U_i\}_{i=1}^n$ are
indeterminates.

Define functions
\begin{eqnarray*}
A\colon \Gen\times\Gen\longrightarrow \Z
&{\text{and}}&
M\colon \Gen\times\Gen\longrightarrow \Z
\end{eqnarray*}
as follows.
Given $\x,\y\in\Gen$, let $$A(\x,\y)=\sum_{i=1}^n
(X_i(\phi)-O_i(\phi)),$$ where $\phi\in\pi_2(\x,\y)$ is any Whitney disk
from $\x$ to $\y$, and $X_i(\phi)$ resp. $O_i(\phi)$ is the algebraic
intersection number of $\phi$ with the submanifold $\{O_i\}\times
\Sym^{g+n-2}(\Sigma)$
resp. $\{X_i\}\times
\Sym^{g+n-2}(\Sigma)$. Also, let
$$M(\x,\y)=\Mas(\phi)-2\sum_{i=1}^n O_i(\phi),$$ where $\Mas(\phi)$
denotes the Maslov index of $\phi$; see~\cite{LipshitzCyl} for an
explicit formula in terms of data on the Heegaard diagram.  Both
$A(\x,\y)$ and $M(\x,\y)$ are independent of the choice of $\phi$ in
their definition.

Moreover, there are functions $A\colon \Gen\longrightarrow \Z$ and
$M\colon \Gen\longrightarrow \Z$ both of which are uniquely specified
to overall translation by the formulas
\begin{eqnarray}
  \label{eq:RelToAbsolute}
A(\x)-A(\y)=A(\x,\y)
&{\text{and}}&
M(\x)-M(\y)=M(\x,\y)
\end{eqnarray}

Let $\UnparModFlow(\phi)$ denote the moduli space of
pseudo-holomorphic disks representing the homotopy class $\phi$,
divided out by the action of $\R$. We can construct an orientation
system $\epsilon$ which gives for each $\phi\in\pi_2(\x,\y)$ with
$\Mas(\phi)=1$ an integer $\#\UnparModFlow(\phi)$, thought of as a suitable
signed count of elements in $\ModFlow(\phi)$.  The counts are gotten by
thinking of $\epsilon$ as coming from an orientation on zero-dimensional
moduli spaces, which are compatible with orientations of the one-dimensional moduli
spaces as follows. Consider the boundaries of all one-dimensional moduli
spaces connecting $\x$ to $\w$ (with $\x\neq \w$). The ends of this moduli space
are of the form
$$\coprod_{
\Big\{
\begin{tiny}
\begin{array}{l}
\y\in \Ta\cap\Tb \\
\phi_1\in\pi_2(\x,\y) \\
\phi_2\in\pi_2(\y,\w)
\end{array}
\end{tiny}
\Big|
\Mas(\phi_1)=\Mas(\phi_2)=1\Big\}}\UnparModFlow(\phi_1)\times\UnparModFlow(\phi_2),$$
where the ends corresponding to $\phi_1, \phi_2$.  For the orientation
of $\ModFlow(\x,\w)$, we have that the induced orientation on the ends
is given by the product orientation on the zero-dimensional moduli
spaces.  Moreover, in the special case where $\x=\w$, there are
additional boundary components, coming from degenerations of the
constant flowline from $\x$ to $\x$ with a holomorphic disks whose
boundary lie entirely within $\Ta$ or $\Tb$.  Our orientation is
chosen so that the degenerations with boundary in $\Ta$ are oriented
with $+1$, while the ones in $\Tb$ appear with orientation $-1$, see~\cite{MOST}.

Let $\CFKm(\oLink)$ be the free module over $\Z[U_1,...,U_n]$
generated by $\Gen$. This module inherits a (relative) bigrading from
the functions $M$ and $A$ above, with the additional convention that
multiplication by $U_i$ drops the Maslov grading by two, and the
Alexander grading by one.

In the case where $\oLink=\oKnot$ has a single component, we define
the differential $$\partial\colon \CFKm(\oKnot)\longrightarrow
\CFKm(\oKnot)$$
by the formula:
\begin{equation}
\label{eq:DefPartial}
\partial (\x)=\sum_{\y\in\Gen}\,
\sum_{\left\{\phi\in\pi_2(\x,\y)\big| 
\begin{tiny}
\begin{array}{l}
\Mas(\phi)=1 \\
X_i(\phi)=0~~~~~~~~~ \forall i=1,...,n
\end{array}
\end{tiny}
\right\}}\!\!
\#\UnparModFlow(\phi)\cm U_1^{O_1(\phi)}\cdots U_n^{O_n(\phi)}\cm \y.  
\end{equation}
It is sometimes convenient
to consider instead the complex $\CFKa(\oKnot)=\CFKm(\oKnot)/(U_1=0)$.
The homology groups
$\HFKm(\oKnot)=H_*(\CFKm(\oKnot))$ and
$\HFKa(\oKnot)=H_*(\CFKa(\oKnot))$ are knot invariants~\cite{Knots},
\cite{RasmussenThesis}, see also~\cite{Links}, \cite{MOS} for the 
case of multiple basepoints, and also~\cite{MOST} for a further discussion
of signs. The gradings induce bigradings
\begin{eqnarray*}
\HFKm(\oKnot)=\bigoplus_{m,s}\HFKm_{m}(\oKnot,s) &{\text{and}}&
\HFKa(\oKnot)=\bigoplus_{m,s}\HFKa_{m}(\oKnot,s).
\end{eqnarray*}

We have defined the bigradings only up to additive constants.  This
indeterminacy can be removed with the following conventions.  Dropping
the condition that all the $X_i(\phi)=0$ in the differential for
$\CFKa$, we obtain another chain complex which retains its Maslov
grading, and whose homology is isomorphic to $\Z$. The convention that
this generator is supported in Maslov grading equal to zero.  Moreover
the Alexander grading can be uniquely pinned down by the requirement
that for each $s$, 
$$\sum_m (-1)^m \rk~\HFKa_{m}(K,s)
=\sum_m (-1)^m \rk~\HFKa_{m}(K,-s).$$

\subsection{Twisted coefficients}
\label{subsec:TwistedCoefficients}

We recall knot Floer homology with twisted coefficients.  Fix an
additional collection of points ${\mathbb P}=\{P_1,...,P_m\}$ on
$\Sigma-\alpha_1-...-\alpha_{g+n-1}-\beta_1-...-\beta_{g+n-1}$ called
a {\em marking}.  We define $P(\phi)$ to be the sum of the local
multiplicities of $\phi$ at the $P_1,...,P_m$.

We can now generalize the earlier construction. 
Set $\RingM=\Z[U_0,...,U_n,t]$, and let $\CFKmT$ be the
free $\RingM$-module generated by $\Gen$,
equipped with the differential
\begin{equation}
        \label{eq:DefTwisted}
\dT(\x)=\sum_{\y\in\Gen}\,
\sum_{\left\{\phi\in\pi_2(\x,\y)\big| 
\begin{tiny}
\begin{array}{l}
\Mas(\phi)=1 \\
X_i(\phi)=0~~~~~~~~~ \forall i=1,...,n
\end{array}
\end{tiny}
\right\}}\!\!
\!\!
\#\UnparModFlow(\phi)\cm t^{P(\phi)}\cm U_1^{O_1(\phi)}\cdots \cm U_n^{O_n(\phi)}\cm \y,
\end{equation}
This generalizes the earlier construction, in the sense that
$$\CFKm(\oKnot)=\CFKmT(\oKnot)/(t-1).$$
The reader should be cautioned, though, that $\CFKm(\oKnot,\Ps)$,
thought of as a module over $\RingM$ is not a knot invariant. However,
if we allow ourselves to invert $t$, we obtain the following:

\begin{lemma}
  \label{lemma:TwistedCoefficients}
  If ${\mathbb P}$ is a marking, then Equation~\eqref{eq:DefTwisted}
  determines a chain complex. Moreover, we have that
  $H_*(\CFKmT(\oKnot)\otimes_{\Z[t]}{\Z[t,t^{-1}]})\cong
  H_*(\CFKm(\oKnot))\otimes_{\Z} \Z[t^{-1},t]$.
\end{lemma}

\begin{proof}
  The fact that it is a chain complex follows along usual lines~\cite{Links},
  in view of the fact that $P$ is additive under juxtaposition of Whitney disks.
  Let $G_s$ denote the set of elements in $\CFKmT(X)$
  of the form $U_1^{k_1}\cm...\cm U_n^{k_n}\cm \x$ for 
  some sequence of integers $k_1,...,k_n$ and $\x\in\Ta\cap\Tb$,
  satisfying the condition that
  $A(\x)-k_1-...-k_n=s$.
  We define a function
  $\epsilon\colon G_s \times G_s\longrightarrow \Z$ as follows.
  Let $U_1^{a_1}\cm...\cm U_n^{a_n} \cm {\mathbf a}$ and 
  $U_1^{b_1}\cm...\cm U_n^{b_n} \cm {\mathbf b}$ be two elements of $G_s$.
  It is easy to see that there is a unique $\phi\in\pi_2({\mathbf a},{\mathbf b})$ with 
  \begin{eqnarray*}
    X_i(\phi)=0 &{\text{and}}&
    O_i(\phi)=b_i-a_i
  \end{eqnarray*}
  for $i=1,...,n$.
  We then define
  $$\epsilon(U_1^{a_1}\cm...\cm U_n^{a_n} \cm {\mathbf a},
  U_1^{b_1}\cm...\cm U_n^{b_n} \cm {\mathbf b})=P(\phi)$$ for this homotopy class.
  Next, fix some ${\mathbf a}_0\in G_s$.
  It is now straightforward to check that the map
  $$\Phi\colon \CFKmT(X) \longrightarrow \CFKm(X) \otimes
  \Z[t^{-1},t]$$
  defined by $$\Phi(U_1^{k_1}\cm...\cm U_n^{k_n}\cm
  \x)=t^{\epsilon({\mathbf a}_0,U_1^{k_1}\cm...\cm U_n^{k_n}\cm \x)}\cm \x$$
  is an
  isomorphism of chain complexes over $\Z[t^{-1},t]$; i.e.  $\Phi$ is
  an isomorphism of $\Z[t^{-1},t]$-modules, and
  $\Phi(\dT\x)=\partial(\Phi(\x))$.
\end{proof}

\begin{figure}[ht]
\mbox{\vbox{\epsfbox{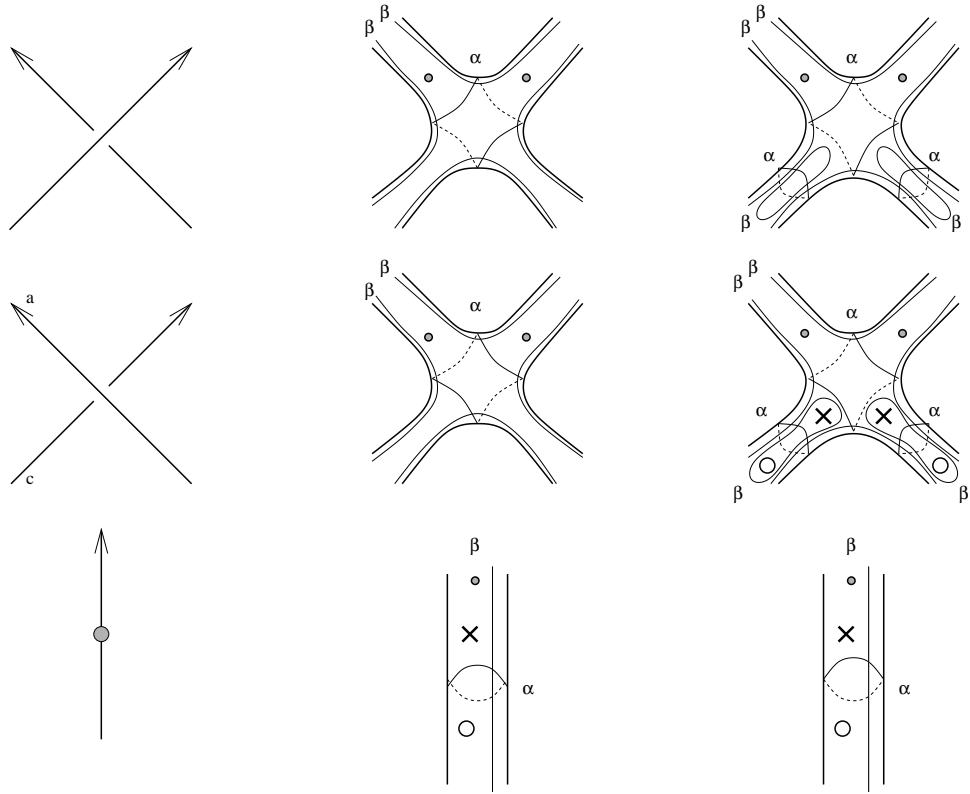}}}
\caption{\label{fig:InitialKnotDiagram}
  {\bf{Heegaard diagram for a decorated projection.}}  The two
  crossings on the left are replaced by a local pictures in the second
  column.  The distinguished edge is replaced by its corresponding
  picture on the bottom. Finally, one of the extra $\beta$ circles
  encircling of the two regions adjacent to the distinguished edge is
  removed. The arc on the left column of the third line represents the
  marked edge, which corresponds to the cylinder in the Heegaard
  surface represented in the left column, decorated with the
  illustrated $\alpha$-circle, and the $O$ point, $X$ point, and
  marking point (represented by the gray circle).  The stabilized
  initial diagram is illustrated in the third column.  }
\end{figure}

For knots in decorated braid position, we can consider the following
Heegaard diagram, which was considered in~\cite{AltKnots}. The
Heegaard surface is a regular neighborhood of the complete
singularization of $K$. All the compact regions in the knot projection
correspond to $\beta$-circles, and the $\alpha$-circles correspond to
crossings.  There is also one additional $\alpha$-circle which serves
as a meridian for the marked edge. We place markings $\Ps$, one at the
top of each edge. We have illustrated this in
Figure~\ref{fig:InitialKnotDiagram}.

Note that in our discussion we have reversed the role of the $\alpha$
and the $\beta$-circles from those in~\cite{AltKnots}. The advantage
is that now the Heegaard surface is endowed with its orientation now
as the outward (rather than inward) boundary of the regular
neighborhood of the singularization. The price is that our formulas
for the Alexander and Maslov gradings associated to states are
slightly different from those in~\cite{AltKnots}, cf.
Figures~\ref{fig:MaslovGradings} and~\ref{fig:AlexanderGradings}
below.

Actually, we will find it useful to have a multiply-pointed Heegaard
diagram, gotten by stabilizing the initial diagram. In the stabilized
diagram, we introduce a new pair of circles for each (unmarked) edge,
one new $\alpha$-circle, which is a meridian, and another new $\beta$-circle,
which is null-homtopic, meeting the corresponding $\alpha$ circle in exactly
two points. The disk bounded by this $\beta$-circle is divided in two
by the $\alpha$-circle, and one of these two regions is marked with
an $X$, the other with an $O$.
We call this diagram the {\em stabilized initial diagram}.

\subsection{Links}

In the case where $\oLink$ is disconnected, some care must be taken in
defining the Heegaard diagram to ensure finiteness of the sum defining
the boundary operator, cf.~\cite{Links}. More precisely, we have the
following:

\begin{defn}
        A {\em periodic domain} is a two-chain $F$ in $\Sigma$
        whose boundary can be expressed as a sum of curves among
        $\alphas$ and $\betas$, and whose multiplicity at $\Os$ and $\Xs$ 
        is zero.
\end{defn}

A Heegaard diagram for an oriented link is said to be {\em admissible}
in the case where all non-trivial periodic domains have both positive
and negative local multiplicities.  Starting with an admissible
Heegaard diagram for a link, we can define $\CFLm$ as above. In this
case, the analogue of $\HFKa$, $\HFLa$, is gotten as follows.  Number
the $\Os$ so that the first $\ell$ of them $\{O_i\}_{i=1}^\ell$ are in
one-to-one correspondence with the different components of the link,
and let $\CFLa=\CFLm/\{U_i=0\}_{i=1}^\ell$. According to~\cite{Links},
$\HFKm(\oKnot)=H_*(\CFKm(\oKnot))$ and
$\HFKa(\oLink)=H_*(\CFLa(\oKnot))$ are link invariants.

Note that link Floer homology from~\cite{Links}, is a $\Z\oplus
\Z^\ell$-graded theory, with one Maslov and $\ell$ Alexander gradings;
the variant we consider here is the $\Z\oplus \Z$-collapsed version,
where the Alexander grading components are added up. We can form the
analogous collapse $\HFKa(\oLink)$, and this is the version which
appears in the skein theory for knots. In fact, a different
construction of the $\Z\oplus\Z$-graded theory appearing in the skein
sequence is given in~\cite{Knots}; the identification of the two
constructions is established
in~\cite[Section~\ref{Links:sec:Relate}]{Links}.

Distinguishing a fixed component of $L$, we could also consider
$H_*(\CFLa(L)/U_1)$, where $U_1$ is the variable corresponding to the
distinguished component.  This group depends on the link only through
the choice of distinguished component.

We can also perform the construction of twisted coefficients as
before, but now the isomorphism from
Lemma~\ref{lemma:TwistedCoefficients} does not hold; indeed, the group
$\HFKmT(\oLink)$ is no longer a link invariant (it depends on the
twisting).

\subsection{Singular links}
In~\cite{SingLink}, the construction of knot Floer homology is
generalized to the case of singular links. This invariant is gotten by
relaxing the constraints on $\Os$ and $\Xs$.  More precisely, for a
link with double-point singularities, the pointed Heegaard diagram
$(\Sigma,\alphas,\betas,\Os,\Xs)$ has the property that
each region $R$ of $\Sigma-\alpha_1-...-\alpha_{g+n-1}={\mathfrak
  A}_1\coprod...\coprod {\mathfrak A}_n$ and
$\Sigma-\beta_1-...-\beta_{g+n-1}={\mathfrak B}_1\coprod...\coprod
{\mathfrak B}_n$
\begin{itemize}
  \item the intersection of $R$ with $\Xs$ contains one or two points
  \item the intersection of $R$ with $\Os$ contains one or two points
  \item the number of points in $R\cap \Xs$ agrees with the number of 
    points in $R\cap \Os$.
\end{itemize}

We can then construct an oriented singular link $S$ in $Y$ (which we
will again assume to be $S^3$ in our applications).  Some pairs of
points in $\Xs$ can be connected by arcs which are disjoint from all
the $\alpha_i$ and $\beta_j$. We call these pairs among the $\Xs$ {\em
  inseparable pairs}. These correspond to the singular points of our
link. The relative Alexander grading $A$ is defined as before, and the
relative Maslov grading is defined by
the formula
$$M(\x,\y)=\Mas(\phi)-2\left(\sum_i O_i(\phi)\right) + 2\left(\sum_j
  XX_j(\phi)\right),$$
where $XX_j(\phi)$ denotes the local
multiplicity of $\phi$ at the inseparable pair $XX_j$. As before, we
can lift the relative Alexander and Maslov gradings to absolute ones,
defined up to overall shifts, 
as in Equation~\eqref{eq:RelToAbsolute}.

We find it convenient also to use a third {\em algebraic grading}
\begin{equation}
  \label{eq:DefAlgebraicGrading}
 N(\x)=M(\x)-2A(\x),
\end{equation}
characterized up to an overall shift by the formula
$$N(\x)-N(\y)=\Mas(\phi)-2 \sum X_i(\phi) +2 \sum XX_j(\phi).$$

With these definitions in place, we can define a chain complexes
$\CFKm(S)$ and $\CFKa(S)$ using a differential as in
Equation~\eqref{eq:DefPartial}. The homology groups of these complexes
depend on the oriented singular link. According to~\cite{SingLink},
\begin{equation}
  \label{eq:EulerSingLink}
  \Delta_S(T)=(T^{-\OneHalf}-T^{\OneHalf})^{\sigma-1} \cm \sum_{d,s}(-1)^d~\rk(\HFKa_d(S,s))\cm T^s,
\end{equation}
where $\sigma$ is the number of singular points in the projection of
$S$, and $\Delta_S$ is the
Alexander polynomial of the singular link, characterized by the skein
relation of Equation~\eqref{eq:SingularSkein}.

The indeterminacty of the $A$-grading can be pinned down by the convention that
$$\sum_{d,s} (-1)^d~\rk\HFKa_d(S,s)=
\sum_{d,s} (-1)^d~\rk\HFKa_d(S,-s).$$ The indeterminacy of the
$M$-grading is pinned down by actually fixing the $N$-grading as
follows. Consider the specialization of the complex to $U_i=1$.
According to \cite{Links}, the resulting complex has the homology of a
torus of dimension $n-1$, which inherits the $N$-grading.  We pin down
the $N$ grading by declaring the bottom-most generator of the torus to
have $N$-grading equal to zero.

A few remarks are in order now, concerning the comparison of this
construction with the conventions adopted in~\cite{SingLink}.  There
are two versions of knot Floer homology for singular links constructed
in~\cite{SingLink}. One of these, ${\mathrm{HFS}}$ is gotten by taking
the above chain complex $\CFm$, identifying the pairs of
indeterminates $U_a$ and $U_b$ which belong to to edges emanating from
the same vertex, setting the indeterminate $U_1$ corresponding to the
first edge to zero, and then taking homology. The second,
${\widetilde{HFS}}$ is gotten by setting each indeterminate $U_i$ to
zero, and then taking homology. The identification used in setting
$U_a=U_b$ has the effect of dividing the Euler characteristic of the
resulting theory by $(1-T)^\sigma$ (which we then choose to
symmetrize); this is how Equation~\eqref{eq:EulerSingLink} in the form
we have stated it follows from results of~\cite{SingLink}.

\section{Calculation for singularized links}
\label{sec:CalcSingLink}

The aim of this section is to calculate a suitably normalized version
of the knot Floer homology groups of a totally singular link in braid
form. 

Let $S$ be a totally singular link in braid form.  We consider
\[
C'(S)=\begin{CD} \CFKmT(S)\otimes_{\RingM} \Big(\bigotimes_{s\in
    \Sing(S)} \RingM @>{t\cm U_a^{(s)}+t \cm U_b^{(s)}-
    U_c^{(s)}-U_{d}^{(s)}}>>\RingM)\Big) ,
  \end{CD}
\]
where $\Sing(S)$ denotes the set of singular points in $S$, the
twisted coefficients for $\CFKmT(S)$ are to be taken with respect to
the special markings illustrated in Figure~\ref{fig:InitialDiagram}
below, and the tensor product is taken (over $\RingM$) with a number
of two-step chain complexes, one for each singular point in the
projection.  As we shall see, these are the chain complexes which
belong at the vertices of the cube of resolutions established in
Theorem~\ref{thm:SpectralSequence}, in the next section.

\begin{figure}[ht]
\mbox{\vbox{\epsfbox{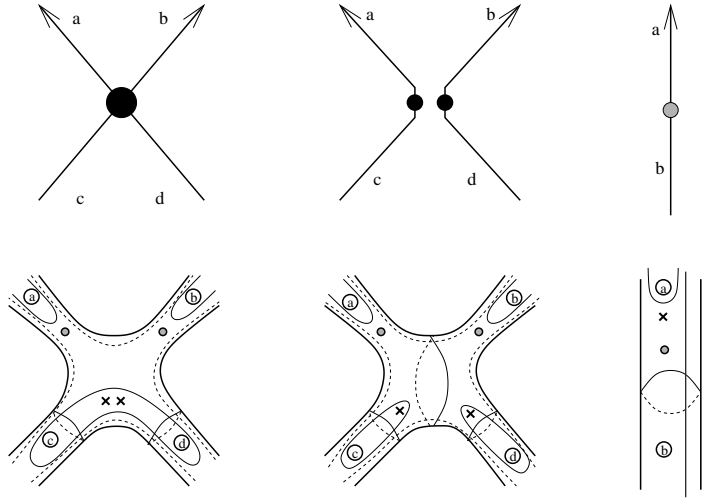}}}
\caption{\label{fig:InitialDiagram} {\bf{Initial diagram for resolved
      links.}} Consider a complete resolution of some knot. We
  construct the Heegaard diagram indicated above (analogous to the
  Heegaard diagram from Figure~\ref{fig:InitialKnotDiagram}.) Note the
  gray dots, representing markings $\Ps$.}
\end{figure}

\begin{theorem}
  \label{thm:CalcSingLink}
  Let $S$ be a singularized link. Then, there are identifications
  \begin{eqnarray}
    \label{eq:Identifications}
    H_*(C'(S))\cong \AlgM(S);
  \end{eqnarray}
  indeed, the homology of $H_*(C'(S))$ is concentrated in a single
  algebraic grading.
\end{theorem}

We prove Theorem~\ref{thm:CalcSingLink} in two steps.  First, we prove
that $H_*(C'(S))$ is all concentrated in a single algebraic grading
(in the sense of Equation~\eqref{eq:DefAlgebraicGrading}). This
follows quickly from an interpretation of generators of this chain
complex in terms of Kauffman states for the projection; indeed, a
version of this is proved in~\cite{SingLink}. The second step
involves a planar Heegaard diagram for the singular link, where the
calculation of $\HFKa(S)$ can be completed readily.

\subsection{Kauffman states}

We state the interpretation of generators in terms of Kauffman states
more explicitly. Fix an oriented diagram for a link, with a
distinguished edge. The projection divides the two-sphere into $n+1$
regions, two of which adjoin the distinguished edge. Let $R(L)$ denote
the remaining regions.

\begin{figure}[ht]
  \mbox{\vbox{\epsfbox{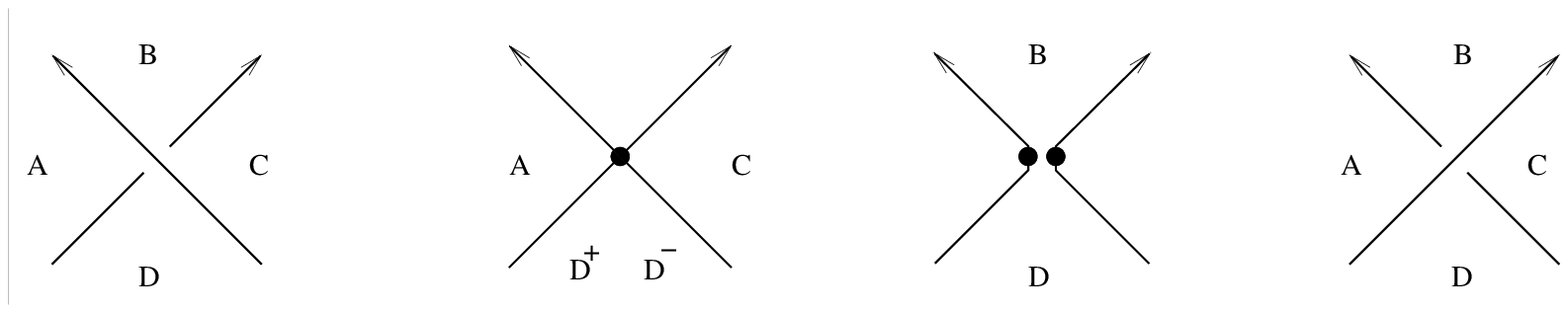}}}
\caption{\label{fig:KauffmanCorners}
  {\bf{Generalized Kauffman corners.}}  We have illustrated here the
  allowed Kauffman corners at a crossing.}
\end{figure}

We define a set, the set of {\em Kauffman corners at $p$}, for each
crossing $p$ of $L$, which depends on whether $p$ is singularized,
resolved, or left alone. If $p$ is on ordinary crossing, the Kauffman
corners $\kappa(p)$ are the four corners of the crossing $p$ $A$, $B$,
$C$, or $D$ (cf.  Figure~\ref{fig:KauffmanCorners}). If $p$ is a
resolved crossing, the Kauffman corners are the two regions $B$ or
$D$. If $p$ is a singularized crossing, the Kauffman corners take
values in $A$, $D+$, $D^-$, and $C$, where both $D^+$ and $D^-$ belong
to the bottom corner $D$. A {\em generalized Kauffman state} for a
singularized link is a map $x$ which associates to each crossing $p\in
c(\Link)$ one of the allowed Kauffman corners $\kappa(p)$, with the
constraint that in each allowed region in $R(L)$, there is a unique
Kauffman corner in the image of $x$, compare~\cite{Kauffman}. We
denote the set of generalized Kauffman states by $\Kauff$.

\begin{figure}[ht]
\mbox{\vbox{\epsfbox{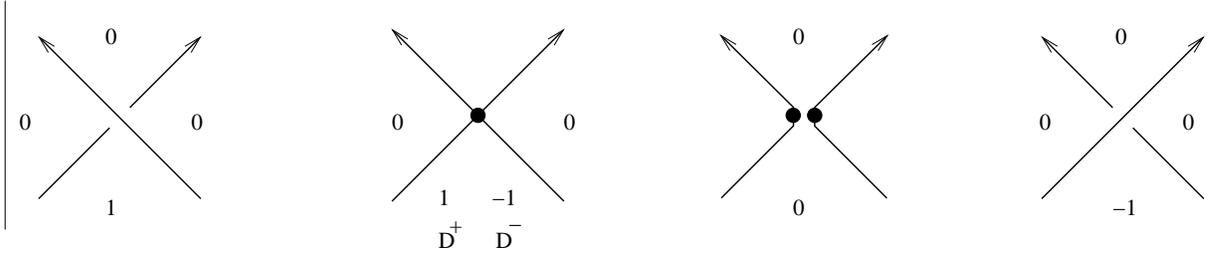}}}
\caption{\label{fig:MaslovGradings}
{\bf{Maslov gradings.}}
This picture illustrates the Maslov gradings for the various Kauffman
corners.}
\end{figure}

\begin{figure}[ht]
\mbox{\vbox{\epsfbox{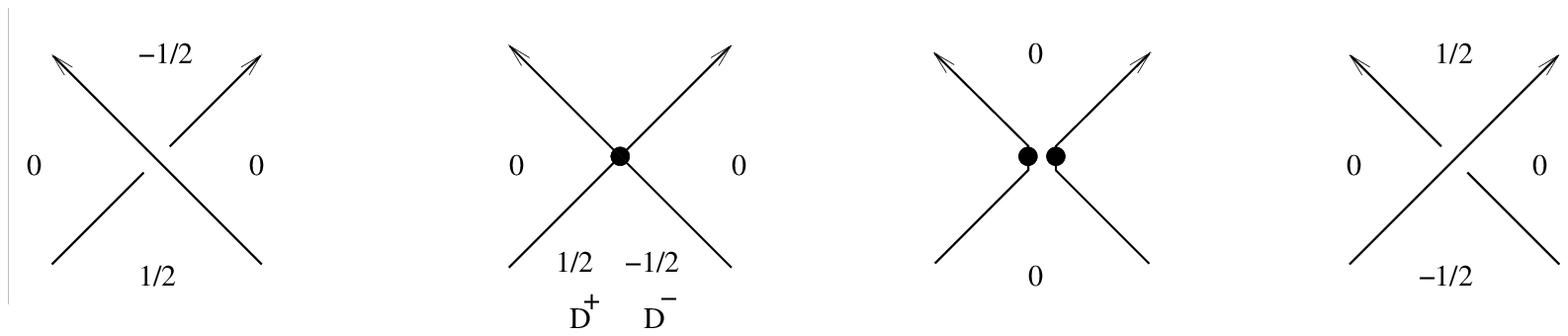}}}
\caption{\label{fig:AlexanderGradings}
  {\bf{Alexander gradings.}}  This picture illustrates the Alexander
  gradings for the various Kauffman corners.}
\end{figure}

Each Kauffman corner has a local Maslov grading $M_p$, which is
assigned as indicated in Figure~\ref{fig:MaslovGradings}.  Similarly,
each Kauffman corner has a local Alexander grading $S_p$, which is
illustrated in Figure~\ref{fig:AlexanderGradings}. The Maslov resp.
Alexander grading of a Kauffman state is given as a sum of local
Maslov resp. Alexander gradings.

We construct a Heegaard diagram for the singular link $S$, gotten by
modifying the Heegaard diagram for the knot projection used in~\cite{AltKnots}
and recalled in Subsection~\ref{subsec:TwistedCoefficients}.

Specifically, the Heegaard surface $\Sigma$ is the boundary of a
regular neighborhood of the link. We associate to each connected,
compact region in the complement of the projection a $\beta$-circle,
and at each crossing $p$ we introduce one additional $\beta$-circle
$\beta_p$, and two $\alpha$-circles $\alpha_c$ and $\alpha_d$. The
circles $\alpha_c$ and $\alpha_d$ are meridians of two of the
in-coming arcs, while $\beta_p$ is a null-homotopic circle meeting
only $\alpha_c$ and $\alpha_d$, in two points apiece. The disk bounded
by $\beta_p$ is divided into three regions by $\alpha_c$ and
$\alpha_d$, one of which meets both $\alpha_c$ and $\alpha_d$ -- which
we mark with a pair of marked points $X$ -- one of which meets only
$\alpha_c$ -- which we mark with $O_c$ -- and the final of which meets
only $\alpha_d$ -- which we mark with $O_d$.

Moreover, to each resolved point, we associate a local picture, again
with two $\alpha$-circles $\alpha_c$ and $\alpha_d$ which are
meridians for the in-coming edges, and now one additional
$\alpha$-circle which bounds a disk which is disjoint from the
singular knot. We introduce also two null-homotopic $\beta$-circles,
$\beta_a$ and $\beta_b$, each of which meets $\alpha_c$ and $\alpha_d$
respectively in exactly two points. Moreover, $\alpha_c$ resp.
$\alpha_d$ divides the disk bounded by $\beta_a$ resp. $\alpha_d$ into
two regions, one of which cointains a basepoits $O_c$ resp. $O_d$, and
the other of which contains a marked point $X$. To the marked edge we
associate an annular region as before.

We also include a marking set ${\mathbb P}$ with one point in each region
on the top of the Heegaard surface corresponding to an edge in
the projection.

This Heegaard diagram is illustrated in
Figure~\ref{fig:InitialDiagram}.

We call this the {\em initial diagram} for the singular link. Note
that when the projection is disconnected, this is not an admissible
diagram in the sense of~\cite{HolDisk}. In particular, it is not clear
that the sums going into the definition of the differential on
$\CFKm(\oKnot,\Ps)$, Equation~\eqref{eq:DefTwisted}, are positive.
This can be remedied by using the marking, and then 
passing to a Novikov completion, as we shall
see presently.
To this end, it is useful to have the following:

\begin{defn}
  Let $(\Sigma,\alphas,\betas,\ws,\zs)$ be a pointed Heegaard diagram
  for a possibly singular knot or link, equipped with a marking $\Ps$.
  We say that $\Ps$ is an {\em admissible marking} if for each integer
  $n$ and $m$, there are only finitely many periodic domains $\pi$
  whose local multiplicities are bounded below by $m$ and for which
  $P(\pi)\leq n$.
\end{defn}

\begin{lemma}
        \label{lemma:Admissibility}
        Let $S$ be a totally singular link.  For the initial diagram,
        the group of periodic domains is the free $\Z$-module
        generated by closed paths in the projection which do not pass
        through the distinguished vertex, i.e. if the projection has
        $n$ components, then we obtain a free Abelian group of rank
        $n-1$.  Moreover, if $K_1,...K_{n-1}$ are these components,
        and $\pi_1,...,\pi_{n-1}$ are its corresponding periodic
        domains, then for $a_i\in\Z$
        \begin{equation}
          \label{eq:WeightFormula}
          P(\sum_{i} a_i\cm \pi_i)=\sum a_i |K_i|,
        \end{equation}
        where $|K_i|$
        denotes the weight of the collection of vertices on $K_i$.
        It follows that the marking on the initial diagram of a singular
        link is an admissible marking.
\end{lemma}

\begin{proof}
  By inspecting Figure~\ref{fig:InitialDiagram}, it is easy to
  construct a periodic domain $\pi_i$ for each knot component $K_i$
  which does not pass through the distinguished vertex.  All the local
  multiplicities of $\pi_i$ are one or zero. They vanish in the
  regions around edges corresponding to components $K_j$ with $j\neq
  i$, and they have local multiplicity $1$ in regions near the edges
  corresponding to the component $K_i$.
  Equation~\eqref{eq:WeightFormula} now follows readily.
  Admissibility of the marking follows as well: any periodic domain
  $\pi$ can be written in the form $\sum a_i \pi_i$. The lower bound
  on $a_i$ the local multiplicities of $\pi$ translates into the bound
  $a_i\geq m$. Finiteness then follows from
  Equation~\eqref{eq:WeightFormula}.
\end{proof}

We can now define the differential on $\CFKmT(S)[[t]]$ as before using the
formula
\[
\dT(\x)=\sum_{\y\in\Gen}\,
\sum_{\{\phi\in\pi_2(\x,\y) | \Mas(\phi)=1, X_i(\phi)=0~~~~~ \forall i=1,...,n\}}\!\!
\#\UnparModFlow(\phi)\cm t^{P(\phi)}\cm U_1^{O_1(\phi)}\cdots U_n^{O_n(\phi)}\cm \y,
\]
only now thinking of the right-hand-side as a formal power series in
$t$.  The fact that for each $t$-power, the sum is finite follows
readily 
Lemma~\ref{lemma:Admissibility}: we are using here that the marking is
admissible. Specifically, the set of all homotopy
classes $\phi\in\pi_2(\x,\y)$ which contribute the same power
$U_1^{k_1}...U_n^{k_n}$ has the form
$$\phi_0+a_1\pi_1+...+a_{n-1}\pi_{n-1},$$
where $\pi_1,...,\pi_{n-1}$
is a basis for the space of periodic domains.  This homotopy class has
a non-zero moduli space only if the corresponding domain has all
non-negative coefficients, which places a lower bound on each of the
coefficients $a_i$. It follows at once that there are only finitely
many terms which contribute for each $t$-power.

Consider next the chain complex
$$C'(S)=\begin{CD} \CFKmT(S)\otimes_{\RingM} \Big(\bigotimes_{s\in
    \Sing(S)} \RingM @>{t\cm U_a^{(s)}+t \cm U_b^{(s)}-
    U_c^{(s)}-U_{d}^{(s)}}>>\RingM)\Big) ,
  \end{CD}$$
using the initial diagram. 

The following is a suitable adaptation of
\cite[Theorem~\ref{SingLink:thm:MoreStates}]{SingLink}.

\begin{prop}
\label{prop:KauffmanStates}
The homology $H_*(C'(S)/U_0)$ is a finitely generated, free
$\Z[[t]]$-module, generated by generalized Kauffman states for the
projection.  All the homology is concentrated in a fixed algebraic
grading, and the Alexander grading is given by the Alexander grading
of the corresponding Kauffman state. Similarly, the homology
$H_*(C'(S))$ is all concentrated in a fixed algebraic grading.
\end{prop}

To establish the above proposition, we employ the following elementary principle:

\begin{lemma}
  \label{lemma:TakeQuotients}
  Let $R=\Z[U_1,...,U_n][[t]]$, and fix an element
  $a=U_1-t\cm \xi$, where $\xi\in R$. Then, 
  the homology of the complex 
  $$
  \begin{CD}
    R@>{a}>>R
  \end{CD}
  $$
  is naturally isomorphic to $R'=\Z[U_2,...,U_n][[t]]$, 
  and the differential (multiplication by $a$) has no kernel.
\end{lemma}

\begin{proof}
  Clearly, $a$ is not a zero-divisor, and hence the homology is
  supported in degree zero.   We find an element $b\in \Z[U_2,...,U_n][[t]]$
  with the property that 
  $$R/(a)\cong R/(U_1-b)\cong \Z[U_2,...,U_n][[t]].$$
  To this end, we are looking for a sequence $\eta_i\in \Z[U_1,...,U_n][[t]]$ with the 
  property that
  \begin{equation}
    \label{eq:WhatWeWant}
    \left(1-\sum_{i=1}^{\infty} \eta_i\cm t^i\right)\cm a = U_1-b.
  \end{equation}
  This is constructed inductively. Indeed, we inductively find 
  \[\{\eta_i\}_{i=1}^m \subset \Z[U_1,...U_n][[t]],\hskip1cm \{b_i\}_{i=1}^m\subset  \Z[U_2,...,U_n],
  \hskip1cm {\text{and}} \hskip1cm R_{m+1}\in \Z[U_1,...,U_n][[t]]\]
  so that
  \begin{equation}
    \label{eq:Induction}
    \left(1-\sum_{i=1}^{m}\eta_i\cm t^i\right) a = U_1-\left(\sum_{i=1}^m b_i t^i\right) + t^{m+1} R_{m+1}.
  \end{equation}
  Having found the first $m$ of the $\eta_i$ and $b_i$ as above, we proceed as follows.
  Clearly, there is some $\eta_{m+1}\in \Z[U_1,..,U_n][[t]]$ with
  $$R_{m+1}(U_1,...,U_n)-R_{m+1}(t\cm \xi,U_2,...,U_n)=\eta_{m+1}\cm a.$$
  Moreover, $R_{m+1}(t\xi,U_2,...,U_n)$ has the following form 
  $$R_{m+1}(t\xi,U_2,...,U_n)=b_{m+1}(U_2,...,U_n) + t \cm
  R_{m+2}(U_1,...,U_n,t),$$
  with $b_{m+1}\in \Z[U_2,...,U_n]$. (The
  induction starts at $m=0$, with $\eta_0=b_0=0$, and $R_{1}=\xi$.)

  Equation~\eqref{eq:WhatWeWant} follows from
  Equation~\eqref{eq:Induction}, taking $m$ to infinity, and letting
  $b=\sum_{i=1} b_i \cm t^i$. Note that the element
  $1-\sum_{i=1}^\infty \eta_i t^i$ is clearly a unit in R, and hence
  the ideal generated by $a$ coincides with the ideal generated by
  $U_1-b$, whose quotient in turn is clearly isomorphic to $R'$.
\end{proof}

We will also need to better understand the generators $\Gen(X)$ of the
initial Heegaard diagram for $X$.  Specifically, note that there is a
map
$$F_X\colon \Gen(X) \longrightarrow \Kauff, $$
which carries a generator to its underlying Kauffman state.
Its fiber consists of $2^m$ points, where here 
here $m$ denotes  the number of singular crossings
plus twice  the number of resolved crossings.

\begin{lemma}
  \label{lemma:Alexander}
  Let $X$ be a singular link.  Given a Kauffman state $k\in
  \Kauff(X)$, there is a unique $\x\in F_X^{-1}(k)$ with maximal
  Alexander grading, denoted $\x(k)$. We claim that there are for any
  state $k\in \Kauff(X)$, the Alexander and the relative Maslov
  gradings of $\x(k)$ are given by the state sum formulas for the
  Alexander and Maslov gradings of $k$
\end{lemma}

\begin{proof}
  This is a straightforward adaptation of
  \cite[Theorem~\ref{SingLink:thm:MoreStates}]{SingLink}.
\end{proof}

\begin{proof}[of Proposition~\ref{prop:KauffmanStates}]
  Let $\Gen=\Ta\cap\Tb$ be the generators for the above chain complex
  $C'(S)/U_0$. Consider the map $F$ from generators $\Gen$ to Kauffman
  states $\Kauff$.  In the special case where the projection of $S$ is
  disconnected, the set of Kauffman states is empty, and hence so is
  the set of generators.  In particular, $H_*(C'(S)/U_0)=0$, and all
  the other statements in the proposition follow at once.
  
  Consider next the case where the projection is connected.  Placing
  basepoints in every region except those disks bounded by
  $\beta$-curves (and also the component containing $U_0$), as
  pictured in Figure~\ref{fig:PointedInitialDiagram} we obtain a new
  chain complex $C''$ whose homology we can calculate.  Specifically,
  we can think of $C''$ as the graded object associated to some
  filtration of $C(S)$. More specifically, if $Q$ denotes the new set
  of auxilliary basepoints, we define the filtration difference
  between generators $\x$ and $\y$ to be the multiplicity of $Q$
  inside the $\phi\in\pi_2(\x,\y)$ with $X(\phi)=0$. (This class
  $\phi$ is uniquely determined by $\x$ and $\y$ since, by hypothesis,
  $S$ corresponds to a singular knot, not a link.)  The differentials
  in the associated graded object now count those homotopy classes
  which are disjoint from $Q$. Indeed, it is easy to see that
  differentials in $C''$ connect two generators only if they
  correspond to the same Kauffman state.  Now, if $\Gen(k)$ is the set
  of generators corresponding to a Kauffman state $k$, then in fact
  the corresponding summand has the form
  $$ 
  \begin{CD}
    N=\bigotimes_{s\in\Sing(S)} \Big(\RingM @>{U_x^{(s)}}>> \RingM\Big)
  \end{CD}
  $$
  where here $x=c$ or $d$.
  We wish to tensor this with
  $$ 
  M=
  \begin{CD}
    \bigotimes_{s\in\Sing(S)} \Big(\RingM @>{t\cm U_a^{(s)}+t \cm U_b^{(s)}-
      U_c^{(s)}- U_d^{(s)}}>> \RingM\Big).
  \end{CD}
  $$

\begin{figure}[ht]
\mbox{\vbox{\epsfbox{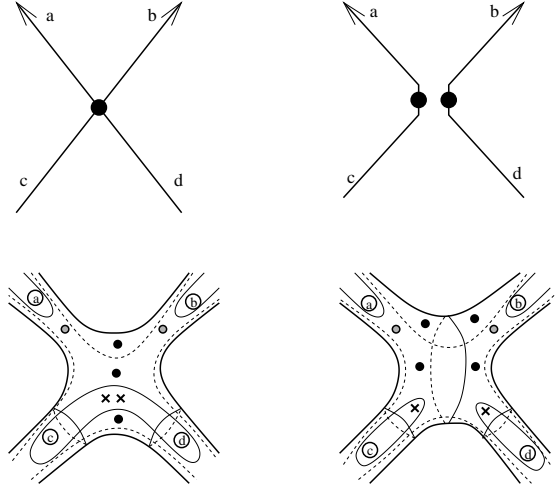}}}
\caption{\label{fig:PointedInitialDiagram}
  {\bf{Pointed initial diagram.}} The black dots represent additional basepoints to be
  added to the initial diagram for the singularized link.}
\end{figure}

  We claim that 
  $$H_*(N\otimes M)\cong \Z[[t]],$$
  as a $\Z[[t]]$-module This is seen
  by repeated applications of Lemma~\ref{lemma:TakeQuotients}.  We
  begin by demonstrating that $H_*(N)\cong \Z[U_1,...,U_m][[t]]$,
  where the variables $\{U_1,...,U_m\}$ are renumbered so that
  each corresponds to a different vertex in the singular link. The point here
  is that of the various $U_x^{(s)}$ we divide out by in $N$, no two
  ever coincides (since they all point into different vertices), thus
  the sequence $\{U_x^{(s)}\}_{s\in X(S)}$ forms a regular sequence,
  i.e. $U_x^s$ is a zero-non-divisor in $\RingM/\{U_x^{v}\}_{v<s}$. It
  follows at once that the homology of $N$ is concentrated in degree
  zero, and it is isomorphic to
  $$\Z[U_1,...,U_{2n}][[t]]/\{U^{(s)}_x\}_{s\in \Sing(S)}\cong
  \Z[U_1,...,U_m][[t]].$$ Letting
  $$ 
  M_{\sigma}=
  \begin{CD}
    \bigotimes_{s\in\Sing(S), s\leq \sigma} \Big(\RingM @>{t\cm U_a^{(s)}+t \cm U_b^{(s)}-
      U_c^{(s)}- U_d^{(s)}}>> \RingM\Big),
  \end{CD}
  $$
  so that 
  $H_*(N\otimes M_0)=H_*(N)$ and 
  $H_*(N\otimes M_\ell)=H_*(N\otimes M)$
  we prove by induction on $\sigma$ that
  $$H_*\left(M\otimes \frac{\Z[U_1,...,U_{2n}][[t]]} {\{t\cm
      U^{(s)}_a+t\cm U^{(s)}_b- U_c^{(s)}- U_d^{(s)}\}}_{s\in
      \Sing(S), s\leq \sigma \}}\right)\cong
  \Z[U_{\sigma+1},...,U_{2n}][[t]]
  $$
  The reason for this is that, once again, $t\cm
  U_a^{(\sigma+1)}+t\cm U_b^{(\sigma+1)}-
  U_c^{(\sigma+1)}-U_d^{(\sigma+1)}$ is a non-zero divisor in
  $\Z[U_{\sigma+1},...,U_{\ell}]$ -- exactly one of $U_c^{(\sigma+1)}$
  or $U_d^{(\sigma+1)}$ vanishes, as it was divided out in $H_*(N)$
  already; the other $U_a^{(\sigma+1)}$ or $U_b^{(\sigma+1)}$ clearly
  does not appear amongst $U_1,...,U_\sigma$. We can now apply
  Lemma~\ref{lemma:TakeQuotients} to verify the inductive step.
  
  The above remarks give us that the homology of $C''/U_0$ (thought of
  as a bigraded module over $\Z[[t]]$) is isomorphic as to the free
  $\Z[[t]]$-module generated by Kauffman states. Indeed, for a
  Kauffman state $k$, the corresponding homology generator is
  represented by the generator with maximal Alexander grading among
  all $\x\in F_X^{-1}(k)$.  Thus, applying Lemma~\ref{lemma:Alexander}
  and glancing at the state sum formula, we see that the Maslov
  grading of $H_*(C''/U_0)$ is twice its Alexander grading, and in
  particular, and hence all the higher differentials in the spectral
  sequence $H_*(C''/U_0)\Rightarrow H_*(C'/U_0)$ (which \em{a priori}
  preserve Alexander grading and drop Maslov grading by one) must
  vanish.

  For the remark concerning $H_*(C'(S))$, observe that there is a long
  exact sequence
  $$
  \begin{CD}
    ...@>>>H_*(C'(S))@>{U_0}>> H_*(C'(S)) @>>> H_*(C'(S)/U_0)@>>>...
  \end{CD}
  $$
  Since multiplication by $U_0$ preserves {\em algebraic} grading, it follows
  that the above sequence splits into direct summands
  $$
  \begin{CD}
    0@>>>H_a(C'(S))@>{U_0}>> H_a(C'(S)) @>>> H_a(C'(S)/U_0)@>>>0,
  \end{CD}
  $$
  where $H_a(C'(S))$ refers to the homology in algebraic grading
  given by $a$, and also $H_a(C'(S)/U_0)=0$ except for $a=a_0$.  Since
  the Alexander grading on $H_a(C'(S))$ is bounded above, it follows
  that if multiplication by $U_0$ (which lowers Alexander grading) is
  an automorphism on $H_a(C'(S))$, then in fact $H_a(C'(S))=0$.
\end{proof}

\subsection{The planar diagram}
\label{subsec:PlanarDiagram}

Proposition~\ref{prop:KauffmanStates} gives an explicit calculation of
$\HFKa(S)$ as a bigraded module over $\Z[[t]]$.  To calculate its
structure as a module over the ring $\RingM$, we find it convenient to
pass to the {\em planar diagram} associated to a singular link,
defined presently, cf.
also~\cite[Section~\ref{SingLink:sec:Planar}]{SingLink}.

As the name suggests, the underlying Heegaard surface in this case is
the sphere.

Each singular point $p$ corresponds to a pair of circles $\beta_p$ and
$\alpha_p$.  The circle $\beta_p$ meets also the $\alpha$-circles
corresponding to the two singular points pointed out of from $p$, in
two points apiece. Indeed, the disk bounding $\beta_p$ is divided in
two by $\alpha_p$, and one of these regions contains the points $O_a$
and $O_b$ corresponding to the out-going edges from our singular
point, while the other region is marked with $XX$. Similarly,
$\beta_p$ divides the disk around $\alpha_p$ into two regions, one of
which we have already encountered, and it is marked by $XX$, and the
other contains the points $O_c$ and $O_d$ corresponding to the
in-coming edges. These regions are also marked with two marking points in
the $\Ps$.

Each smoothed point corresponds to a pair of $\alpha$-circles,
$\alpha_a$ and $\alpha_b$ and a pair of $\beta$-circles $\beta_c$ and
$\beta_d$, where $\alpha_a$ and $\beta_c$ meet in two points, and
$\alpha_b$ and $\beta_d$ meet in two points. Letting $D_a$, $D_b$,
$D_c$, and $D_d$ be the disks bounded by $\alpha_a$, $\alpha_b$,
$\beta_c$, and $\beta_d$ respectively, we have that $D_a\cap D_c$ and
$D_b\cap D_d$ are marked by $X$, and $D_a-D_a\cap D_c$, $D_b-D_b\cap
D_d$, $D_c-D_a\cap D_c$, and $D_d-D_b\cap D_d$ resp. are marked by
$O_a$, $O_b$, $O_c$, and $O_d$ respectively.

\begin{figure}[ht]
\mbox{\vbox{\epsfbox{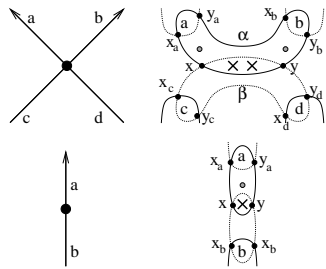}}}
\caption{\label{fig:LocalPicture}
{\bf{Diagram near the singularized or resolved crossing.}}
        Intersection points here are labeled.}
\end{figure}

At the distinguished edge, we have a $\alpha$-circle corresponding to the 
out-going edge, and an $\beta$-circle corresponding to the in-coming edge.
These two circles are in fact disjoint, and the complementary region is
marked with a single $X$, as illustrated in Figure~\ref{fig:LocalPicture}.

We include also markings contained inside the disks bounded by the
$\alpha$-circles.

\begin{figure}[ht]
\mbox{\vbox{\epsfbox{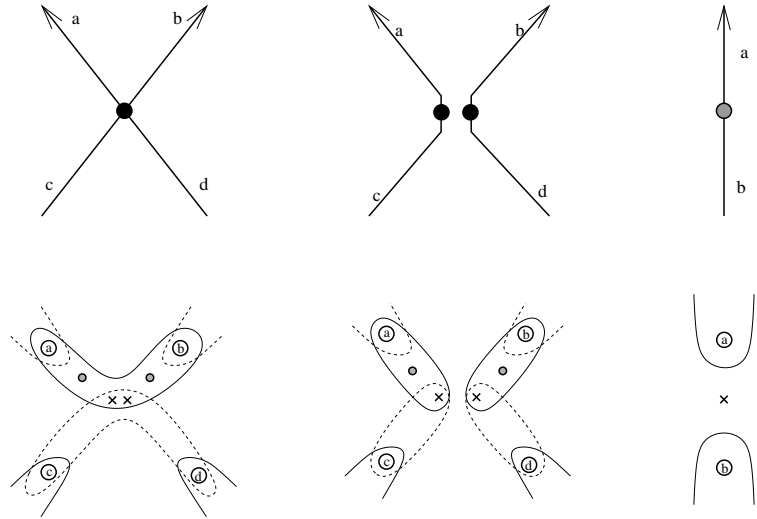}}}
\caption{\label{fig:PlanarDiagram}
  {\bf{Planar diagram for singularized links.}}  Given part of a
  singular knot projection illustrated on the top row, we
  construct the following planar diagram illustrated on the bottom row. 
  We have illustrated from the top row: a singularized point, a smoothed point,
  and the distinguished edge.}
\end{figure}

This Heegaard diagram is illustrated in Figure~\ref{fig:PlanarDiagram}.

For this diagram, we can form the same sum as in
Equation~\eqref{eq:DefTwisted}. The present diagram satisfies a stronger
finiteness condition: each non-trivial periodic domain has both positive and
negative local multiplicities. As in~\cite{HolDisk}, this ensures that
the differential from Equation~\eqref{eq:DefTwisted} consists of
finitely many terms, and hence it can be done over the base ring
$\RingM$. Let $\CFKmTP(S)$ denote the resulting chain complex.

\begin{lemma}
  \label{lemma:HandleslideInvariance}
        There is an isomorphism
        $H_*(\CFKmTP(S)[[t]])\cong H_*(\CFKmTI(S)[[t]])$,    
        where $\CFKmTI(S)$ is the chain complex for a singularized knot
        defined using the initial diagram, and $\CFKmTP(S)[[t]]$
        is the chain complex using the planar diagram, and completed at $t$.
\end{lemma}

\begin{proof}
  We describe a sequence of handleslides and destabilizations going
  from the initial diagram to the planar diagram. Consider first the
  case where there are no smoothed edges. Order the edges in the
  diagram, starting at the distinguished edge, and proceeding
  according to the orientation of the knot. This then assigns numbers
  to the two outgoing edges $a$ and $b$ and the two incoming edges $c$
  and $d$.  At a given singular point, we take the outgoing edge which
  is assigned the highest number, and handleslide its corresponding
  meridian ($\alpha$-circle) over the other three meridians
  corresponding to the three other edges of the singular point. In
  fact, we do this in the following order: we start with the very last
  edge in the diagram, and work our way backwards. It is easy to see
  that what we are left with can be destabilized to go back to the
  planar diagram.
  
  We can now follow the handleslides and destabilizations with
  isomorphisms among knot Floer homology groups,
  following~\cite{HolDisk}. Handleslide invariance is established via
  a map which counts pseudo-holomorphic triangles. The key point now
  is that the handleslides we indicated in the above algorithm never
  cross the marked points $\Ps$.
  
  More specifically, let $(\Sigma,\alphas,\betas,\ws,\zs)$ and
  $(\Sigma,\alphas',\betas,\ws,\zs)$ be two pointed Heegaard diagrams,
  where $\alphas$ and $\alphas'$ differ by handleslides; and suppose
  that $\Ps$ is a marking which is admissible for both diagrams, so
  that the marking $\Ps$ is disjoint from the support of the
  handleslides. This is the condition which ensures that
  $$H_*(\CFKmT(\Ta,{\mathbb T}_{\alpha'}))\cong H_*(T^{g+n})\otimes
  \RingM;$$
  so that it has a canonical (up to sign) top-dimensional
  generator $\Theta$.  Following
  \cite[Section~\ref{HolDisk:sec:HandleSlides}]{HolDisk},  we can
  define the handleslide map
  $${\underline F}\colon \CFKmT(\Sigma,\alphas,\betas,\ws,\zs)[[t]] \longrightarrow
  \CFKmT(\Sigma,\alphas',\betas,\ws,\zs)[[t]].$$
  by
  $${\underline F}(\x)=\sum_{\y\in{\mathbb T}_{\alpha'}\cap\Tb}
  \sum_{\{\psi\in\pi_2(\x,\Theta,\y)\big|\Mas(\psi)=0, X_i(\psi)=0\}}
  \#\ModFlow(\psi) \cm t^{P(\psi)}\cm U_1^{O_1(\psi)}\cdot...\cdot
  U_n^{O_n(\psi)}\cm \y,$$
  where here as usual $\pi_2(\x,\Theta,\y)$
  denotes the space of homotopy classes of Whitney triangles
  connecting $\x\in\Ta\cap\Tb$, $\Theta\in\Ta\cap{\mathbb
    T}_{\alpha'}$, and $\y\in{\mathbb T}_{\alpha'}\cap\Tb$.  Indeed,
  exactly as in
  \cite[Section~\ref{HolDisk:sec:HandleSlides}]{HolDisk}, it follows
  that this map induces an isomorphism in homology.
  quasi-isomorphism.
  
  Invariance under destabilizations follows immediately from the usual
  proof, cf.~\cite[Section~\ref{HolDisk:sec:Stabilization}]{HolDisk}.

  The case where there are smoothings follows similarly. In fact, the
  only point here is to pay extra attention to the ordering of the
  edges: we take the ordering of the edges starting at the
  distinguished edge, and then proceeding according to the orientation
  of the resolved singular knot.
\end{proof}

For the planar diagram of a singular knot, we organize the
intersection points $\Ta\cap\Tb$ as follows. At each singularized
point $p$, the intersection points of $\beta_p$ with $\alpha_q$ with
$\alpha_q\neq \alpha_p$ can be partitioned into two pairs, each of which
is indexed by one of the outgoing edges $a$ or $b$ from $p$. We label
the points $x_p$, $y_p$, $x_a$, $y_a$, $x_b$, and $y_b$.  We say that
$x_a$ and $y_a$ are assocated to the edge $a$, and $x_b$ and $y_b$ are
associated to the edge $b$. Similarly, at each of the two vertices $v$
for a smoothing, we have a pair of circles $\beta_v$ and $\beta_v$,
which intersect at two points which we label $x$ and $y$. Similarly,
the intersection points of $\beta_v$ with $\alpha_q$ with $\alpha_q\neq
\alpha_v$, we have a pair of intersection points which we label $x_a$
and $y_a$, which are associated to the edge $a$.  These labeling
conventions are illustrated in Figure~\ref{fig:LocalPicture}.

Given an generator $\x$, we consider the graph which is the union over
all $x_i\in \x$ of the edges associated to $x_i$.

\begin{defn}
  \label{def:Type}
  Let $S$ be a singular link. A {\em coherent cycle} is an closed,
  connected cycle embedded in the graph underlying the singular link
  which can be oriented coherently with the orientation of $S$, and
  which does not include the distinguished edge.  A {\em coherent
    multi-cycle} is a disjoint union of coherent cycles, or possibly
  the empty set.
\end{defn}

\begin{lemma}
        \label{lemma:AssociatedGraph}
  The graph corresponding to any generator is a coherent multi-cycle.
  Morover, for a fixed multi-cycle, there are $2^n$ different generators
  which map to it. 
\end{lemma}

\begin{proof}
  This is a straightforward consequence of the combinatorics of our
  Heegaard diagram.
\end{proof}

There is a canonical generator $\x_0$, which is the product of the
$x_p\in\alpha_p\cap\beta_p$ (as in Figure~\ref{fig:LocalPicture}) for
all $p$.  The corresponding graph is empty.

\begin{figure}[ht]
\mbox{\vbox{\epsfbox{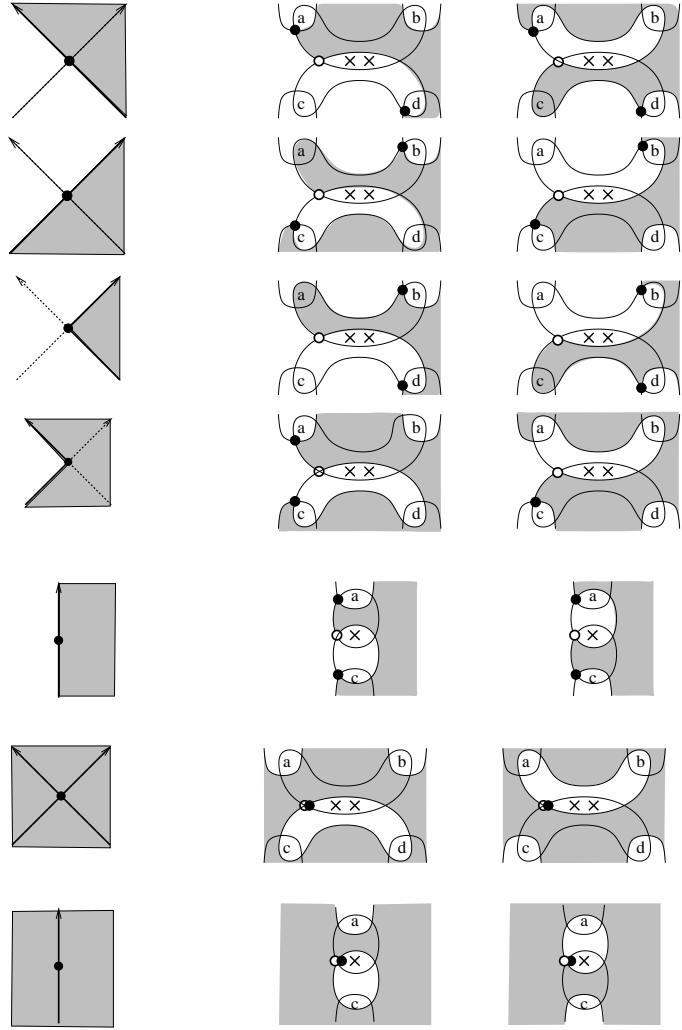}}}
\caption{\label{fig:PlanarDomains} {\bf{Instructions for constructing
      the two possible cyclic domains.}} The left hand column consists
  of the possible local pictures of a cycle $Z$ in the knot
  projection, where the dark edges signify edges connecting two
  vertices contained inside the the region bounded by $Z$ (while the
  dotted edges do not); moreover, the interior of the cycle is shaded
  gray. The second column signifies the corresponding two-chain for
  $\phi_1$, while the third indicates the corresponding two-chain for
  $\phi_2$.}
\end{figure}

Consider the case where $Z$ consists of a single component, and
suppose also that $S$ has connected underlying graph.  Consider the
corresponding generator $\x(Z)$ which contains only points of the form
$x_p$, $x_a$, or $x_c$ (i.e. it contains none of the points of the
form $y_p$, $y_a$, or $y_c$ at any $p$), and whose associated
multi-cycle is $Z$.  We claim that there are exactly two positive
domains $\phi_1$ and $\phi_2$ in $\pi_2(\x(Z),\x_0)$ with
$X_i(\phi)=0$.  The precise recipe for constructing both $\x(Z)$ and
the domains $\phi_1$ and $\phi_2$ are illustrated in
Figure~\ref{fig:PlanarDomains}.  To see that there are only two, we
proceed as follows. Any domain in $\pi_2(\x(Z),\x_0)$ can be written
as $\phi_1+\sum_{s\in\Sing(S)} k_s \cm P_s$, where $k_s$ are some
integers, and $P_s$ is the periodic domain at the singular point $s$.
By looking at the local picture near each singular point, we see that
if this domain is positive, then $|k_s|\leq 1$. In fact, if $k_s\neq
0$, and if $e$ is an edge connecting two vertices inside $Z$ from
$s_1$ to $s_2$, then $k_{s_1}=k_{s_2}$.  It now follows easily that
$\phi_1$ and $\phi_2$ are the only positive domains.

The fact that $n_X(\phi_1)=n_X(\phi_2)=0$ follows from the fact that
our planar diagram is in braid form, with distinguished edge on the
left. It is easy to see that $\Mas(\phi_i)=1$. (In fact, more will be
proved in Lemma~\ref{lemma:HolomorphicRepresentatives}.)

\begin{prop}
  \label{prop:MinimalElement}
  The element $\x_0$ is a cycle.
\end{prop}

\begin{proof}
  Clearly, $\x_0$ has minimal algebraic degree among all generators
  which correspond to the empty graph.

  In the above discussion, we have shown that if $Z$ is a connected
  cycle, then the generator $\x(Z)$ with minimal algebraic grading among
  all generators corresponding to $Z$ has algebraic degree one greater than
  that of $\x_0$. Iterating this, it follows readily that
  if $Z$ has $n$ components, then any of its corresponding generators
  has algebraic grading at least $n$ greater than that of $\x_0$.
\end{proof}

To complete the proof of Theorem~\ref{thm:CalcSingLink}, then, it
suffices to verify that the boundaries of the various chains with
algebraic grading one greater than that of $\x_0$ give the stated
relations.  Relations of the form 
\begin{equation}
\label{eq:PointRelations}
t^2\cm U^{(p)}_a\cm U^{(p)}_b 
=
U^{(p)}_c U^{(p)}_d
\end{equation}
arise as the boundaries of the various generators
corresponding to the empty graph, which replace exactly one component
of $\x_0$ by the corresponding $y_p$. The remaining relations come
from looking at the boundaries of $\x(Z)$, with the help of the
following:

\begin{lemma}
  \label{lemma:HolomorphicRepresentatives}
  Consider a Heegaard diagram, and suppose that there are two generators $\x$ and $\y$
  for its Heegaard Floer homology which can be connected by a homology class of disks
  $\phi\in\pi_2(\x,\y)$ whose associated two-chain $\cald(\phi)$ has the following properties:
  \begin{itemize}
  \item  $\cald(\phi)$ is a planar region, with one boundary consisting of a $2n$-gon,
    where alternate between arcs in $\{\alpha_i\}_{i=1}^n$ and $\{\beta_i\}_{i=1}^n$
    (thus the corners alternate between components of $\x$ and components of $\y$)
  \item all other boundaries of $\phi$ are bounded by circles $\{\alpha_j\}_{j=n+1}^{n+m}$
  \item there are arcs in the $\{\beta_j\}_{j=n+1}^{n+m}$ which
    connect between various of the $\{\alpha_i\}_{i=1}^{n+m}$ with the
    condition that, $\cald(\phi)-\beta_{n+1}-...-\beta_{n+m}$ has a
    single connected component; in fact, if we construct the graph
    whose vertices are arcs in $\alpha_i$ and whose edges correspond
    to $\beta_j$ which connect them, we obtain a graph with a single
    closed cycle, corresponding to the $2n$-gon on the boundary of
    $\cald(\phi)$.   (It is easy to see that all the componets of $\x$ in the interior $\alpha_i$ (with $i>n$)
    coincide with the corresponding components of $\y$.)
  \end{itemize}
  In this case, $\Mas(\phi)=1$, and for any choice of almost-complex
  structure, $\#\UnparModFlow(\phi)=\pm 1$.
\end{lemma}

\begin{proof}
  We can embed $\phi$ into a Heegaard diagram for $S^3$ with the
  following properties. There are only three generators of
  $\Ta\cap\Tb$, one of which is $\x$, the other $\y$, and the third is
  denoted $\w$, and the only positive Whitney disk supported away from the basepoint $O$
  is $\phi$. Since $\HFa(S^3)$ is one-dimensional, it follows at once that $\phi$ has 
  a unique holomorphic representative (for any almost-complex structure).
  
  The Heegaard diagram is constructed by starting with a Heegaard
  diagram for $S^3$, thought of as plumbing of a sequence of $n$
  spheres, the first of which has square $-2$, the next $n-2$ have
  square zero, and the last square $-1$.  This gives us the case $\phi_0$ where
  $m=0$. 
  
  To pass to the more general case, we stabilize the diagram as
  follows.  Suppose $\beta_{j}$ is some arc which connects
  $\alpha_{j}$ $j>n$ to some $\alpha_i$ with $i\leq n$.  We then
  attach a one-handle to the Heegaard surface constructed so far in
  such a manner that one of the feet is supported in the region marked
  by the basepoint $O$, and the other is attached inside the support
  of $\phi_0$. We then add also the curve $\alpha_{j}$, which is a
  core of our one-handle, and complete $\beta_{j}$ to be a closed
  curve which runs through our handle. After adding all the $\alpha_j$
  with distance one from the boundary, we proceed to those with
  distance two. Now, the one-handle is added so that one foot is
  inside the region bounded by the $\alpha$-circle with smaller
  distance.  Proceeding inductively, we construct our desired Heegaard
  diagram for $S^3$, cf. Figure~\ref{fig:PolygonalRegion}.

\begin{figure}[ht]
\mbox{\vbox{\epsfbox{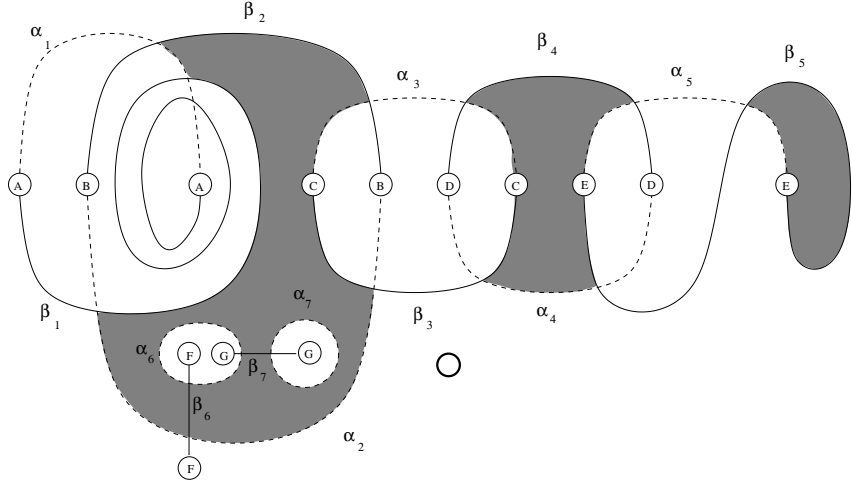}}}
\caption{\label{fig:PolygonalRegion}
  {\bf{A Heegaard diagram for $S^3$.}}  The illustrated Heegaard
  diagram for $S^3$ has the property that $\CFa(S^3)$ has only three
  generators.  One-handles are labelled with letters, the marked point
  is indicated by an unmarked hollow circle.  There is exactly one
  non-negative flowline which does not cross the marked point, and it
  is shaded.  The domain here is a $10$-gon with two disks removed.}
\end{figure}

\end{proof}

In view of Lemma~\ref{lemma:HolomorphicRepresentatives} (and its
analogue, replacing roles of $\alpha$- and $\beta$-circles), taking
differentials of minimal algebraic chains in $\x(Z)$ (where $Z$ is a
coherent cycle), we obtain relations in the homology of the form
\begin{equation}
  \label{eq:IndeterminateSignGenRelZ}
t^{|W_Z|} \cm \prod_{e\in \Out(W_Z)\cap \In(W^c_Z)} U_e
 = \pm \prod_{e\in \In(W_Z)\cap \Out(W^c_Z)} U_e
\end{equation}
Equation~\eqref{eq:GenRel}, where $W_Z$ is the set of vertices
contained inside the multi-cycle $Z$. (More precisely: the complement
of the cycle $Z$ consists of two open sets; one of these regions $D$
contains the distinguished edge, and the other does not. By $W_Z$ we
mean all vertices which in the complement of $D$.)

The sign is determined by the following 

\begin{lemma}
  \label{lemma:Signs}
  Letting $\x_0$ be the minimal intersection point from
  Proposition~\ref{prop:MinimalElement}, and let $\x(Z)$ 
  be the generator with minimal algebraic degree associated
  to a fixed (one-component) cycle $Z$ for the graph. Let
  $\phi_1,\phi_2\in\pi_2(\y,\x_0)$ be the two Maslov index one
  homotopy classes. Then, with respect to the usual sign conventions,
  $\#\UnparModFlow(\phi_1)=-\#\UnparModFlow(\phi_2)$.
\end{lemma}

\begin{proof}
  According to Lemma~\ref{lemma:HolomorphicRepresentatives},
  $|\#\UnparModFlow(\phi_1)|=|\#\UnparModFlow(\phi_2)|=1$. We exclude
  the possibility that $\#\UnparModFlow(\phi_1)=
  -\#\UnparModFlow(\phi_2)$.  Let $C$ denote the chain complex
  associated to the planar diagram, now setting all the $U_i=1$ and
  $t=1$. This homology group is invariant under isotopies of the
  diagram which cross the $O_i$.  Moreover, $\x_0$ remains a cycle in
  in $C$, and the existence of $\y$ as above ensures that the homology
  class represented by $\x_0$ is non-trivial, but has two-torsion.
  But if we isotope the planar diagram (crossing the $O_i$) so that at
  each $\alpha_\sigma$ meets only its corresponding $\beta_\sigma$
  (and none of the other $\beta$-circles), then we are left only with
  intersection points corresponding to the empty graph, and in the
  resulting chain complex, all differentials cancel in pairs. In
  particular, there is no torsion in the homology.
\end{proof}

Our aim is to show that the general case follows from these special
cases.

Lemma~\ref{lemma:Signs} determines the sign in
Equation~\eqref{eq:IndeterminateSignGenRelZ}. We wish to show that all
the other relations in Equation~\eqref{eq:GenRel} follow from
relations the corresponding relations for sets $W=W_Z$, where $Z$
ranges over all coherent cycles $Z$ in the projection, together with
the relations from Equation~\eqref{eq:PointRelations}.  More
precisely, given a coherent cycle $Z$, let $W_Z$ denote the set of
vertices contained in the interior of $Z$. Let $\AlgM'(S)$ denote the
algebra $\Z[U_0,...,U_{2n}]$ divided out by the relations
\begin{equation}
  \label{eq:GenRelC}
t^{|W_Z|} \cm \prod_{e\in \Out(W_Z)\cap \In(W^c_Z)} U_e
 = \prod_{e\in \In(W_Z)\cap \Out(W^c_Z)} U_e
\end{equation}
and 
\begin{equation}
t^2 \cm U_a^{(s)}\cm  U_b^{(s)} = U_c^{(s)}\cm  U_d^{(s)}. 
\end{equation}
at each singular point $s$.
There is an obvious projection map
$$q\colon \AlgM'(S) \longrightarrow \AlgM(S).$$
Indeed, we have the following:

\begin{lemma}
  The map $q$ is an isomorphism.
\end{lemma}

\begin{proof}
  We prove that each relation of the form Equation~\eqref{eq:GenRel}
  follows from corresponding relations, where we take only those sets
  $W'$ which are connected, and have connected complement.  Moreover,
  they follow from corresponding relations where there are no vertices
  with only in-coming or only out-going edges. But the relations
  corresponding to sets with this form are precisely the relations
  from Equation~\eqref{eq:GenRelC}.
\end{proof}

Supposing that $S$ has connected projection, we have shown that that
the homology $\HFKa(S)$ in the Maslov degree given by $\x_0$ coincides
with the $\RingM$-module (which, incidentally, is also a ring) $\AlgM(S)$
described in the introduction.  In the case where the graph underying
$S$ is connected, it is easy to see that $\AlgM(S)$ is non-trivial, as
it contains the element corresponding to $1$; since
Proposition~\ref{prop:KauffmanStates} ensures that all of $\HFKa(S)$
is concentrated in one algebraic grading, we conclude that $\HFKa(S)\cong
\AlgM(S)$.

In the other case where the graph underlying $S$ is
disconnected, $\AlgM(S)$ is trivial, but so is $\HFKa(S)$, as can be
readily seen from Proposition~\ref{prop:KauffmanStates}.
 
\vskip.2cm 
\begin{proof}[Of  Theorem~\ref{thm:CalcSingLink}.]
The above discussion gives the identification 
$$H_*(\CFKmT(S))\cong \AlgM(S).$$
\end{proof}

It is perhaps worth noting that Theorem~\ref{thm:CalcSingLink} gives
an algebraic calculation of the invariant~$\widetilde{\mathrm{HFS}}$
of~\cite{SingLink} for connected, planar singular knots.
Specifically, 
$${\widetilde{HFS}_*}(S)\otimes \Z[t{-1},t]]=
{\mathrm{Tor}}_*^{\Z[U_0,...,U_{2n},t^{-1},t]]}\left(\AlgM(S),\Z[U_0,...,U_{2n},t^{-1},t]]/(U_0,...,U_{2n})\right).$$

\section{Skein exact sequences, and a cube of resolutions}
\label{sec:ExactSeq}

Our aim now is to demonstrate exact sequences relating the knot Floer
homology of a knot, its singularization, and its smoothing.  We begin
with a the following variant of Theorem~\ref{thm:SkeinExactSequenceIntro}:

\begin{theorem}
  \label{thm:SkeinExactSequence}
  Let $\Kp$ be an oriented knot given equipped with a distinguished
  positive crossing $p$, $\Smooth$ denote its (oriented) smoothing at
  $p$, and $\Singularize$ denote its singularization at $p$.  Then,
  for a suitable choice of twisting, we have the following exact
  sequences
  $$
  \begin{tiny}
  \begin{CD}
    ...@>>>\HFKm(\Kp)\otimes \Z[t] @>>>
    H_*(\frac{\CFKmT(\Singularize)}{(t\cm U_a + t\cm U_b - U_c - U_d)}) @>>>
    H_*(\CFLmT(\Smooth))
    @>>>...
  \end{CD}
  \end{tiny}
  $$
  and 
  $$
  \begin{tiny}
  \begin{CD}
    ...@>>>\HFKa(\Kp)\otimes \Z[t] @>>>
    H_*(\frac{\CFKaT(\Singularize)}{(t\cm U_a + t\cm U_b - U_c - U_d)}) @>>>
    H_*(\frac{\CFLmT(\Smooth)}{U_a})@>>>...
  \end{CD}
  \end{tiny}
  $$
  Letting $\Km$ denote the corresponding knot where the
  positive crossing is changed to negative,  we have the exact sequences:
  $$
  \begin{tiny}
  \begin{CD}
    ...@>>>\HFKm(\Km)\otimes \Z[t]@>>> \HFKmT(\Smooth)
    @>>> H_*(\frac{\CFKmT(\Singularize)}{(t\cm U_a + t\cm U_b - U_c - U_d)})
    @>>>...
  \end{CD}
  \end{tiny}
  $$
  and
  $$
  \begin{tiny}
  \begin{CD}
    ...@>>>\HFKa(\Km)\otimes \Z[t]@>>> H_*(\frac{\CFLmT(\Smooth)}{U_a})
    @>>> H_*(\frac{\CFKaT(\Singularize)}{(t\cm U_a + t\cm U_b - U_c - U_d)})
    @>>>...
  \end{CD}
  \end{tiny}
  $$
\end{theorem}

Theorem~\ref{thm:SkeinExactSequence} is stated with twisted
coefficients. Indeed, for the singular link, one is to understand
coefficients twisted using the special markings illustrated
in Figure~\ref{fig:InitialDiagram}. Specializing to $t=1$, we
obtain the following:

\begin{cor}
  \label{cor:SkeinExactSequence}
  With $\Kp$, $\Km$, $\Singularize$, and $\Smooth$ as above, we have
  long exact sequences
  $$
  \begin{CD}
    ...@>>>\HFKm(\Kp) @>>>
    H_*(\frac{\CKm(\Singularize)}{(U_a+U_b- U_c-U_d)}) @>>>
    \HFLm(\Smooth) @>>> ... \\
    ...@>>>\HFKa(\Kp)@>>>
    H_*(\frac{\CFKaT(\Singularize)}{(U_a+U_b-U_c-U_d)}) @>>> 
    H_*(\frac{\CFLm(\Smooth)}{U_a})@>>>... \\
    ...@>>>\HFKm(\Km)@>>> \HFKm(\Smooth)
    @>>> H_*(\frac{\CFKmT(\Singularize)}{(U_a+U_b-U_c-U_d)})
    @>>>... \\
    ...@>>>\HFKa(\Km)@>>> H_*(\frac{\CFKm(\Smooth)}{U_a})
    @>>> H_*(\frac{\CFKaT(\Singularize)}{(U_a+U_b-U_c-U_d)})
    @>>>...
  \end{CD}
  $$
\end{cor}

\begin{figure}[ht]
\mbox{\vbox{\epsfbox{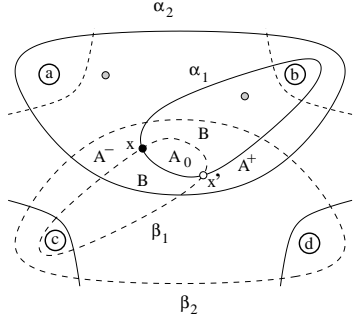}}}
\caption{\label{fig:ExactSequence}
{\bf{Exact triangle.}} Markings near a singular point used in
in the proof of the exact triangle with singular links.}
\end{figure}

Theorem~\ref{thm:SkeinExactSequence} is proved by inspecting a
suitable Heegaard diagram, pictured in Figure~\ref{fig:ExactSequence}.
After proving Theorem~\ref{thm:SkeinExactSequence} and its corollary,
we turn our attention to a stronger form of the theorem which gives
the hypercube of resolutions alluded to in the introduction.  Note
that the version stated in Theorem~\ref{thm:SkeinExactSequenceIntro}
involves a statement about absolue Alexander and Maslov degrees. We
return to this detail in Subsection~\ref{subsec:AbsoluteMaslovAlexander}.

We find it convenient to give a proof of
Theorem~\ref{thm:SkeinExactSequence} in terms of a slight modification
of the planar Heegaard diagram from
Subsection~\ref{subsec:PlanarDiagram}.  Specifically, draw a Heegaard
diagram near a singular point for $\Singularize$ as shown in
Figure~\ref{fig:ExactSequence}. In that picture, we have distinguished
circles $\alpha_1$ and $\beta_1$ which meet in two points $x$ and
$x'$.  If we drop the pair of circles $\alpha_1$ and $\beta_1$, and
then choose $A^-$ and $A^0$ as points in $\Xs$ (and disregard the
points $B$), we obtain a marked Heegaard diagram for $\Singularize$.
This diagram is marked with $\Os=\{O_1,...,O_n\}$, and
$\Xs=\{A^-,A^0\}\cup \Xs_0$, where $\Xs_0=\{X_1,...,X_{n-2}\}$.
Alternatively, if we leave in $\alpha_1$ and $\beta_1$, and use
$\Xs=\Xs_0\cup \{A^-,A^0\}$, we obtain a Heegaard diagram for
$\Km$.  Leaving in $\alpha_1$ and $\beta_1$, and using
$\Xs=\Xs_0\cup \{A^0\cup A^+\}$, we obtain a Heegaard diagram for the
knot with positive krossing $\Kp$.  Finally, using $\Xs$ as the
union of $\Xs_0$ and the two regions marked by $B$, we obtain a
Heegaard diagram for the smoothing $\Smooth$ of the crossing.

Clearly, $\CFKmT(\Km)$ has a subcomplex $X$ consisting of configurations
which contain the intersection point $x$, and a quotient complex $Y$. Thus, 
$\CFKmT(\Km)$ can be thought of as the mapping cone of the map
$$\Phi_{B}\colon Y \longrightarrow X$$
gotten by counting  Maslov index one flowlines
which contain exactly one of the regions marked by $B$, i.e.
$$\Phi_{B}(\x)=\sum_{\y\in{\mathfrak S}} \sum_{\left\{\phi\in\pi_2(\x,\y)
  \big| 
\begin{tiny}
\begin{array}{l}
\Mas(\phi)=1, \\
 X_i(\phi)=0~~~~~~~~~~\forall i=1,...,n-2 \\
B_1(\phi)+B_2(\phi)=1 
\end{array}
\end{tiny}
\right\}} \#\UnparModFlow(\phi) \cm t^{P(\phi)}\cm
U_1^{O_1(\phi)}\cm...\cm U_n^{O_n(\phi)}\cm \y. $$
Moreover,
$\CFKmT(\Smooth)$ has $Y$ as a subcomplex, with quotient $X$, and
hence, it can be thought of as the mapping cone of the map
$$\Phi_{A^-}\colon X \longrightarrow Y,$$
defined by counting
flowlines which contain exactly one of the regions marked by $A^0$ or
$A^-$.  Clearly $X$ is isomorphic to $\CFKmT(\Singularize)$.

Similarly, there is a subcomplex $X'$ of $\CFKmT({\mathcal R})$
consisting of configurations which contain the intersection point
$x'$.  This has a quotient complex we denote by $Y'$. Moreover, $\Kp$
has a subcomplex isomorphic to $Y'$ and quotient complex isomorphic to
$X'$. Thus, we can think of $\CFKmT(\Kp)$ as a mapping cone of
$$\Phi_{B'}\colon X'\longrightarrow Y'$$
gotten by counting flowlines with non-zero multiplicity at $B$, 
and $\CFKmT(\Smooth)$ as the mapping cone of
$$\Phi_{A^+}\colon Y'\longrightarrow X',$$
which counts flowlines through exactly one of $A^0$ or $A^+$.

\begin{lemma}
\label{lemma:CalcComposite}
The composite maps $\Phi_{A^+}\circ \Phi_{B'}$ 
and $\Phi_{B}\circ\Phi_{A^-}$ are chain homotopic to multiplication by
$$\pm (t\cm U_a + t\cm U_b - U_c - U_d).$$
\end{lemma}

\begin{figure}[ht]
\mbox{\vbox{\epsfbox{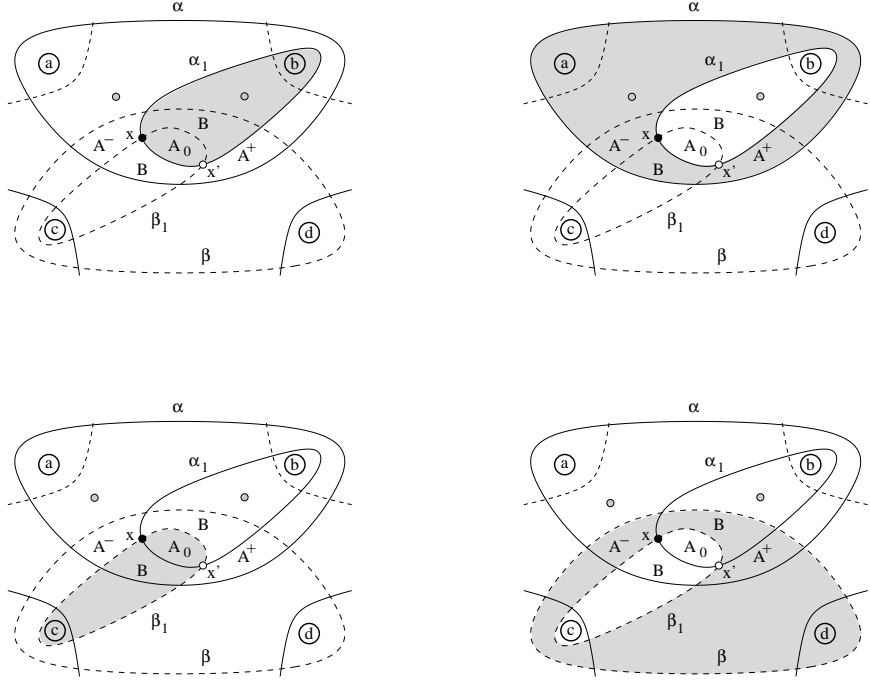}}}
\caption{\label{fig:CalcComposite}
{\bf{Composite of $\Phi_A\circ\Phi_B$.}} This is the geometric content
of Lemma~\ref{lemma:CalcComposite}. The four
types of boundary degeneration with Maslov index equal to two
and which contain one of both $A$ and $B$ are pictured here.}
\end{figure}

\begin{proof}
  Let $\Phi_{A^-B}\colon X \longrightarrow X$ denote the
  map defined by counting flowlines with $\Mas(\phi)=1$ 
  and which represent homotopy classes with the property that
  $A^0(\phi)+A^{-}(\phi)=1$ and also the total multiplicity in $\phi$ of the two
  points labelled by $B$ is equal to one. Considering ends 
  of moduli spaces with Maslov index equal to two, we see that
  $$\partial\circ \Phi_{A^-B}+\Phi_{A^-B}\circ \partial +
  \Phi_{B}\circ \Phi_{A_-} + F=0,$$
  where $F$ is a count of Maslov index
  two boundary degenerations containing one of the points among
  $\{A_0,A_-\}$ and one of the two regions marked by $B$. But there
  are four such homotopy classes, two from
  $\Sigma-\alpha_1-...-\alpha_n$, and two from
  $\Sigma-\beta_1-...-\beta_n$. The former two contribute $t\cm U_a+t
  \cm U_b$ and the latter two $-U_c-U_d$, cf.
  Figure~\ref{fig:CalcComposite}.  Note that the signs work as stated:
  the first two terms come from boundary degenerations with boundary
  on $\Ta$, and the second from those with boundary on $\Tb$.

  The case of $\Phi_{A^+}\circ \Phi_{B'}$ follows similarly.
\end{proof}

\begin{proof}[Of Theorem~\ref{thm:SkeinExactSequence}]
We can now consider the following two filtered complexes, denoted
${\mathcal P}$ and ${\mathcal N}$ respectively:
\begin{eqnarray}
\begin{diagram}
X' &\rTo^{\Phi_{B'}} & Y' \\
\dTo^{t\cm U_a + t\cm U_b - U_c - U_d} & \rdTo^{\Phi_{A^+B}} & \dTo_{\Phi_{A^+}} \\
X' & \rTo_{=} & X'
\end{diagram}
&\text{{\hskip.5cm}and{\hskip.5cm}}&
\begin{diagram}
X &\rTo^{=} & X \\
\dTo^{\Phi_{A^-}} & \rdTo^{\Phi_{A^-B}} & \dTo_{t\cm U_a + t\cm U_b - U_c - U_d} \\
Y & \rTo_{\Phi_{B}} & X
\end{diagram}
\label{eq:defPN}
\end{eqnarray}

The first of these is clearly quasi-isomorphic to $\CFKmT(\Kp)$. 
This can be seen by
considering the filtration induced by the vertical coordinate, where
the complex is expressed as the extension of a subcomplex with trivial
homology by a quotient complex which is $\CFKmT(\Kp)$. 

Alternatively, if
we consider the filtration by the horizontal coordinate, where the
first complex has a subcomplex quasi-isomorphic to $\CFKmT(\Smooth)$, 
with quotient
complex the mapping cone of
$$t\cm U_a + t\cm U_b - U_c - U_d \colon X' \longrightarrow X'.$$

Putting these two observations together, we get a long exact sequence:
\[
  \begin{CD}
    ...@>>> H_*(\CFKmT(\Smooth))   @>>> H_*(\CFKmT(\Kp))  @>>> H_*(\frac{X'}{t\cm U_a+t\cm U_b-U_c-U_d}) @>>>
    \end{CD}
\]

A similar reasoning applied to the second complex gives the following:
\[
  \begin{CD}
    ...@>>>  H_*(\frac{X}{t\cm U_a+t \cm U_b-U_c-U_d}) @>>> H_*(\CFKmT(\Km)) @>>> H_*(\CFKmT(\Smooth)) @>>>
  \end{CD}
\]
Observe now that $X$ is identified with $\CFKmT(\Singularize)$, while
$X'$ is also identified with $\CFKmT(\Singularize)$ (only with a shift
in gradings, cf. also Subsection~\ref{subsec:AbsoluteMaslovAlexander}
below).  The two exact triangles involving $\HFKa(\Kp)$ and
$\HFKa(\Km)$ are gotten by setting $U_a=0$.
\end{proof}

\begin{proof}[Of Corollary~\ref{cor:SkeinExactSequence}.]
  The corollary follows at once by specializing the above discussion
  to $t=1$.
\end{proof}

If $\Knot$ is a (possibly singular) knot, and $p$ is a positive
crossing, whose singularization is denoted $s_p(\Knot)$ and whose
smoothing is denoted $r_p(\Knot)$, then we can define chain maps
\begin{equation}
  \label{eq:DefUnzipDiag}
\Unzip_p\colon \CFKmT(s_p(\Knot)) \otimes_{\RingM}
\RingM/(t\cm U_a^{(p)}+t\cm U_b^{(p)}- U_c^{(p)}-
U_d^{(p)})\longrightarrow \CFKmT(r_p(\Knot)).
\end{equation}
and
\begin{equation}
  \label{eq:DefZipDiag}
  \Zip_p\colon \CFKmT(r_p(\Knot)) \longrightarrow \CFKmT(s_p(\Knot))\otimes_{\RingM}
\RingM/(t\cm U_a^{(p)}+t\cm U_b^{(p)}-U_c^{(p)}-U_d^{(p)}),
\end{equation}
which we call ``unzip'' and ``zip'' homomorphisms respectively. These are
the homomorphisms correponding to the horizontal arrows from 
Diagram~\eqref{eq:defPN}.

Let $\Knot$ be a knot projection, and consider its cube or resolutions
$\bigcup_{I\colon c(\Knot)\longrightarrow \{0,1\}} X_I(\Knot)$.  Note
that for the purposes of the results from this section, it is
irrelevant if the projection has braid form or not.

We can form a graded $\RingM$-module
$V(\Knot)=\bigoplus_{I\colon c(\Knot)\longrightarrow \{0,1\}} V(X_I(\Knot))$,
where here 
$$V(X_I(\Knot))=H_*\left(\CFKmT(X_I(\Knot))\otimes
  \Big(\bigotimes_{s\in \Sing(X_I(\Knot))} \RingM/(t\cm U_a^{(s)}+t\cm
  U_b^{(s)}-U_c^{(s)}-U_d^{(s)})\Big)\right),$$
where here $\Sing$
denotes the set of singular points in the singularized knot
$X_I(\Knot)$.

\begin{theorem}
  \label{thm:SpectralSequence}
  Let $\Knot$ be a knot equipped with a projection.
  There is a spectral sequence whose $E_1$ term is
  $V(\Knot)$, and with $d_1$ differential induced by the zip and unzip
  homomorphisms, which converges to $\HFKmT(\Knot)$.
\end{theorem}

Draw the planar Heegaard diagram for the knot $\Knot$, gotten by putting
together the pictures from Figure~\ref{fig:ExactSequence}. In
particular, at each crossing $p$, we have a pair of $\Xs$-markings
$A^0_p$ and $A^-_p$ or $A^+_p$, depending on the sign of the crossing.
In fact, we drop the distinction from our notation, denoting both
markings simply by $A_p$.  We place pairs of markings $B_p$ as well,
also as in Figure~\ref{fig:ExactSequence}.

There is a natural filtration on the chain complex for $\CFKmT(\Knot)$. 
Specifically, let
${\mathfrak S}$ denote the set of generators $\Ta\cap\Tb$ for the
Heegaard diagram for $\Knot$.
Given a crossing $i$, we say that a generator $\x$ has type
$X_i$ if either:
\begin{itemize}
\item $i$  is a negative crossing and $\x$ contains the distinguished point
$x_i$
\item $i$ is a positive crossing and $\x$ contains the distinguished point
$x_i'$
\end{itemize}
A homology class of Whitney disk $\phi$ is
said to have {\em type $A_i$} if the following conditions are
satisfied:
\begin{itemize}
\item $\phi\in \pi_2(\x,\y)$ where one of $\x$ and $\y$ has type $X_i$ and the other has type $Y_i$
\item $\phi$ has multiplicity at one of the two
points marked $A_i$, and zero at the other
\item $\phi$ has multiplicity zero at both points marked $B_i$.
\end{itemize}
Similarly, $\phi$ is said to have {\em type $B_i$} if
\begin{itemize}
\item $\phi\in\pi_2(\x,\y)$, where one of $\x$ and $\y$ has
type $X_i$ and the other has type $Y_i$
\item $\phi$ has multiplicity zero at both points marked with $A_i$.
\item $\phi$ has multiplicity one at one of the two points marked
with $B_i$ and multiplicity zero at the other.
\end{itemize}
Finally, $\phi$ has {\em type $A_i B_i$} if
\begin{itemize}
\item $\phi\in\pi_2(\x,\y)$ where both $\x$ and $\y$ are of type $X_i$
\item $\phi$ has multiplicity one at one of the two regions marked $A_i$,
and zero at the other
\item $\phi$ has multiplicity one at one of the two regions marked $B_i$,
and zero at the other.
\end{itemize}

We have maps 
\begin{eqnarray*}
\Phi_{A_i}&\colon& X_i \longrightarrow  Y_i \\
\Phi_{B_i}&\colon& Y_i \longrightarrow  X_i 
\end{eqnarray*}
gotten by counting pseudo-holomorphic disks of type $A_i$ resp. $B_i$.
More generally, fix disjoint sets $I,J,K \subset \{1,...,n\}$.
Let $X_I=\cap_{i\in I} X_i$, $Y_J=\cap_{j\in J} Y_j$, $X_K=\cap_{k\in K} X_k$.
We can define a map
$$\Phi_{I J K}\colon 
X_I\cap Y_J\cap X_K \longrightarrow Y_I\cap X_J\cap X_K$$
gotten by counting flows which are 
\begin{itemize}
\item of type $A_i$ for all $i\in I\cap P$
and $B_i$ for all $i\in I \cap N$, 
\item type $B_j$ for all $j\in J\cap P$
and $A_j$ for all $j\in J \cap N$
\item 
type $A_k B_k$ for all $k\in K$,
\item the multiplicity at all the other basepoints $A_i$ and $B_i$ are zero
\end{itemize}
where $N$ are the negative crossings and $P$ are the positive ones.

\begin{lemma}
  \label{lemma:Complex1}
        Fix disjoint sets $I,J,K\subseteq\{1,...,n\}$. Then
\begin{eqnarray*}
\lefteqn{       \sum_{
\begin{tiny}
{\left\{\begin{array}{l}
I_1, J_1, K_1,\\ I_2, J_2, K_2
\end{array} \Big|
\begin{array}{l}
                I_1\coprod I_2\coprod K_1 \coprod K_2 = I\cup K \\
                J_1\coprod J_2\coprod K_1 \coprod K_2 = J\cup K
\end{array}\right\}}\end{tiny}}
                 \Phi_{{I_1},{J_1},{K_1}}\circ  
                \Phi_{{I_2},{J_2},{K_2}}
        } 
   \\
&=&
\left\{
\begin{array}{ll}
  -u^{(k)} & {\text{if $K$ contains a single element $k$,
and $I=J=\emptyset$}} \\
0 & {\text{otherwise.}}
\end{array}
\right.
\end{eqnarray*}
where $u^{(k)}$ is multiplication by the scalar $t\cm U^{(k)}_a + t\cm U^{(k)}_b
- U^{(k)}_c - U^{(k)}_d$.
\end{lemma}

\begin{proof}
  This follows from the usual Lagrangian Floer homology proof that
  $\partial^2=0$: we consider ends of moduli spaces of
  pseudo-holomorphic disks with Maslov index equal to two, and
  consider the ends. Multiplication by the scalars arises from the
  boundary degenerations from Lemma~\ref{lemma:CalcComposite}.  Note
  that boundary degenerations with Maslov index two contain exactly
  one $A_i$ and its corresponding $B_i$; i.e. they are of the form
  counted in that lemma.
\end{proof}

Consider the graded group
$$C=\bigoplus_{I\coprod J\coprod K \coprod L=\Sing(\Knot)} X_I\cap Y_J\cap X_K \cap X_L,$$
whose index set consists of all partitions of the crossings into
four disjoint sets. This index set can be ordered by the convention that
$$ (I_1,J_1,K_1,L_1)\preceq (I_2,J_2,K_2,L_2) $$
if $I_2\subseteq I_1$, 
$J_2\subseteq I_1\cup J_1$,
$K_2\subseteq I_1\cup K_1$.
Consider the map
$$D_{(I_1,J_1,K_1,L_1)\preceq (I_2,J_2,K_2,L_2)}\colon C_{I_1, J_1, K_1, L_1} \longrightarrow 
C_{I_2, J_2, K_2, L_2}$$
defined by the following formula
\begin{equation}
\label{eq:DefD}
D_{(I_1,J_1,K_1,L_1)\preceq (I_2,J_2,K_2,L_1)}
=\left\{\begin{array}{ll}
\Phi_{{J_2\cap I_1}, {L_2\cap J_1}, {L_2\cap I_1}}
& {\text{if $K_1=K_2$}} \\
\pm 1 & {\text{if $K_2=I_1\cup \{n\}$, $J_2=J_1$, $L_2=L_1$}} \\
& {\text{and $n$ is a negative crossing}} \\
\pm 1 & {\text{if $L_2=K_1\cup \{p\}$, $J_2=J_1$, $L_2=L_1$}} \\
& {\text{and $p$ is a positive crossing}} \\
\pm u^{(n)} & {\text{if $L_2=K_1\cup \{n\}$, $J_2=J_1$, $L_2=L_1$}} \\
& {\text{and $n$ is a negative crossing}} \\
\pm u^{(p)} & {\text{if $K_2=I_1\cup \{p\}$, $J_2=J_1$, $L_2=L_1$}} \\
& {\text{and $p$ is a positive crossing}} \\
\end{array}
\right.
\end{equation}
The signs $\pm 1$ are determined as follows. They depend on the generator $\x$,
and also on the pair of gradings $(I_1,J_1,K_1,L_1)\preceq (I_2,J_2,K_2,L_2)$,
where $K_1$ and $K_2$ differ by a single crossing $c$. Specifically,
if we order all the crossings in the diagram, then 
\begin{equation}
  \label{eq:SignChoice}
  \pm 1 = (-1)^{M(\x) + \#\{k\in K_1\cap K_2\big| c<k\}}.
\end{equation}

We can endow $C$ with an endomorphism
$$D=\sum_{(I_1,J_1,K_1,L_1)\preceq (I_2,J_2,K_2,L_2)} D_{(I_1,J_1,K_1,L_1)\preceq (I_2,J_2,K_2,L_2)}.$$

\begin{lemma}
  The module $C$ endowed with the endomorphism $D$ is a
  chain complex, which is quasi-isomorphic to $\CFKmT(\Knot)$.
\end{lemma}

\begin{proof}
  The fact that it is a complex follows from
  Lemma~\ref{lemma:Complex1}. More precisely, we must verify that the
  $(I_2,J_2,K_2,L_1)$ component of $D^2|_{(I_1,J_1,K_1,L_1)}$ vanish.
  In the cases where $K_2=K_1$, this is a direct consequence of
  Lemma~\ref{lemma:Complex1}. In the cases where $|K_1|$ and $|K_1|$
  differ in a single element, this is a consequence of the fact that
  any of the $\Phi$-maps commutes (up to sign) with multiplication by
  $1$ or by $u^{(s)}$. The factor $(-1)^{M(\x)}$ in the choice of
  signs is inserted so that these squares, in fact, anti-commute.
  Finally, we must check $D^2=0$ in cases where $K_1$ and $K_2$ differ
  in two places. In this case, we must consider squares where two of 
  the edges are marked by multiplication by $\pm 1$ or $\pm u^{(s)}$.
  Anticommutativity of the squares, then, is ensured by the factor
  $(-1)^{\#\{k\in K_1\cap K_2\big| c<k\}}$ in the signs.
  
  The quasi-isomorphism with $\CFKmT(\Knot)$ is seen by contracting
  horizontal differentials labelled by $\pm 1$ in
  Equation~\eqref{eq:DefD}. Equivalently, each singular point $i$
  endows $C$ with a filtration as in Equation~\eqref{eq:defPN}.
  Contracting those horizontal differentials labelled by $=$ in
  Equation~\eqref{eq:defPN}, we end up with a chain complex for
  $\CFKmT(\Knot)$.  More precisely, consider the subcomplex of $C$
  corresponding to $(I,J,K,L)$ with some positive crossing $p$ is in
  $K$ or $L$. This is clearly a subcomplex, and it has trivial
  homology.  In the quotient complex (which can be thought of as a
  complex we obtain after contracting the bottom horizontal arrow in
  each ${\mathcal P}$), we have a subcomplex consisting of those
  $(I,J,K,L)$ with the property that each negative crossing $n$ is in
  $J$ or $L$. Clearly, its quotient complex has trivial homology. This
  subcomplex, in turn is identified directly with $\CFKmT(\Knot)$.
\end{proof}

\begin{proof}[Of Theorem~\ref{thm:SpectralSequence}]
  We can think of our chain complex as filtered by an $n$-dimensional
  hypercube by collapsing the vertical filtration. The resulting
  chain complex is
  $$
  \bigoplus_{I\colon c(\Knot)\longrightarrow \{0,1\}}
  \CFKmT(X_I(\Knot)),
  $$
  where 
  $$\CFKmT(X_I(\Knot))=
  \CFKmT(X_I(\Knot))\otimes
  \Big(\bigotimes_{s\in \Sing(X_I(\Knot))} \RingM/(t\cm U_a^{(s)}+ t\cm U_b^{(s)}-
  U_c^{(s)}-U_d^{(s)})\Big),$$
  with edge homomorphisms induced by zip and unzip maps.
  This immediately gives rise to our stated spectral sequence.
\end{proof}

\subsection{Absolute Maslov and Alexander gradings and Theorem~\ref{thm:SkeinExactSequenceIntro}.}
\label{subsec:AbsoluteMaslovAlexander}

The above proof of Theorem~\ref{thm:SkeinExactSequence} verifies the
version stated in the introduction,
Theorem~\ref{thm:SkeinExactSequenceIntro}, up to overall shifts in
Alexander and Maslov gradings of the three terms: it is clear that the
maps (which are induced by $\Phi_B$, $\Phi_{A^\pm}$, and connecting
homomorphisms) all preserve Alexander gradings and Maslov gradings, up
to shifts. It remains to pin down this indeterminacy.

\begin{proof}[Of Theorem~\ref{thm:SkeinExactSequence}]
  We re-examine the diagram Equation~\eqref{eq:defPN}, indicating
  Alexander gradings induced from ${\mathcal P}$ in parentheses:
\begin{eqnarray}
\begin{diagram}
X'(s) &\rTo^{\Phi_{B'}} & Y'(s) \\
\dTo^{ U_a +  U_b - U_c - U_d} & \rdTo^{\Phi_{A^+B}} & \dTo_{\Phi_{A^+}} \\
X'({s-1}) & \rTo_{} & X'({s-1}).
\end{diagram}
\end{eqnarray}

We verify that $\phi_+$ preserves Alexander gradings as follows.
Choose an element $\x'\in X'(s)$. This element inherits Alexander
grading $s$ from the complex for $\Kp$. The map $\phi_+$ is induced by
the projection of the above complex to the mapping cone:
$$U_a+U_b-U_c-U_d\colon X'({s}) \longrightarrow X'({s-1}),$$
which in turn is identified with the mapping cone
$$U_a+U_b-U_c-U_d\colon X(s+1) \longrightarrow X(s).$$
Here $X$ is the complex for the singularized link, endowed with the
Alexander grading inherited from the diagram for $\Kp$, when it is
thought of as generated by those elements which contain the the
distinguished intersection point $x$ from
Figure~\ref{fig:ExactSequence}, rather than $x'$, which was used to
define $X'$. This is the convention compatible with the state sum
formula (and hence also it is the Alexander grading which induces the
Alexander polynomial as specified by the skein sequence of
Equation~\eqref{eq:SingularSkein}).  Thus, the element $\x'$ in
$\Kp$-Alexander bigrading grading $s$ projects
to the element $\x'\in X(s+1)$, which we think of as an element
in Alexander grading $s$ in the mapping cone.
Thus, $\phi_+$ preserves the Alexander grading.

To see that $\delta_+$ preserves Alexander gradings, we proceed as
follows.  Modifying our Heegaard diagram for $\Kp$, by moving the two
basepoints $X$ to the new pair $B$ as in the proof of
Theorem~\ref{thm:SkeinExactSequence}, we obtain a new Heegaard diagram
for $\Smooth$. We must compare the internal Alexander grading of this
diagram with the one it inherits by thinking of its chain complex as a
summand of ${\mathcal P}$. Observe first that generators for this new
complex also give rise to Kauffman states for the resolution, except
for states which contain ``right'' Kauffman corners, labelled by $C$
in Figure~\ref{fig:KauffmanCorners}, and those cancel in pairs by a
suitably chosen rectangle crossing none of the basepoints.  Consider
now a generator $\x$ for the $\Smooth$, which corresponds to a
Kauffman state for the resolution. There are two cases: either the
Kauffman state occupies the top corner, or the bottom corner. If it
occupies a top corner, it is an element of $Y'(s)$, as shown in
Figure~\ref{fig:CorrespondingToResolved}, whose Alexander grading is
determined by a state sum formula. Alternatively, it can be considered
as a generator for the positve crossing, where it again corresponds to
a top Kauffman state, and hence its local contribution is $1/2$ (so
that its Alexander grading is $1/2$ greater than when the state was
considered as a generator for the resolved diagram).  When $\y$
occupies a bottom state (in $Y'(s)$) for the resolved diagram, then we
can find another state (as in
Figure~\ref{fig:CorrespondingToResolved}) in $X'({s+1})$, i.e. which
has one greater Alexander grading, and which in fact represents the
Kauffman corner $D_+$ for the resolved diagram. Thus, the Alexander
grading of $\y$, when thought of as a generator for $\Kp$, is once
again $1/2$ greater than its Alexander grading, when thought of as a
generator for $\Smooth$. This same argument readily shows that
$\psi_+$ preserves Alexander gradings.

\begin{figure}[ht]
\mbox{\vbox{\epsfbox{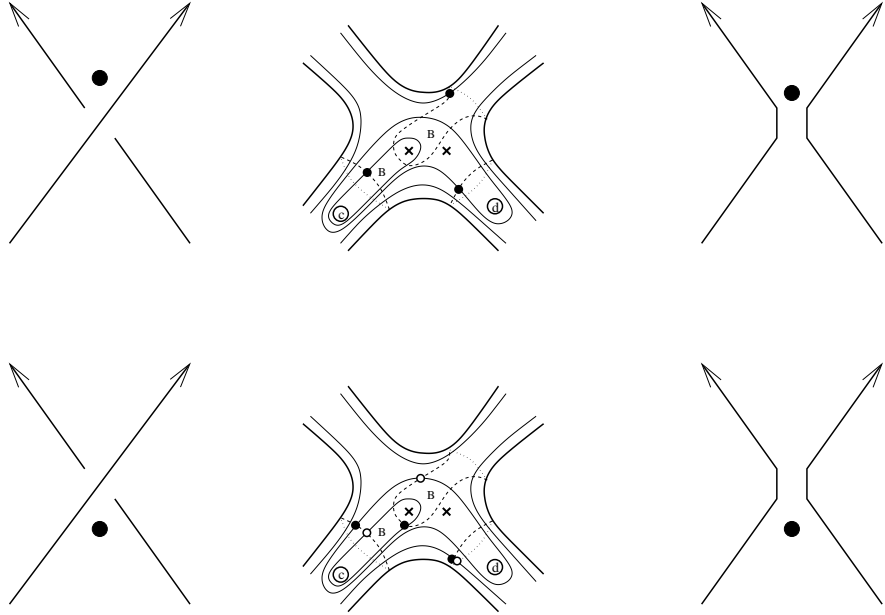}}}
\caption{\label{fig:CorrespondingToResolved} {\bf{Verification of
      Alexander grading shifts.}}  The top drawing demonstrates a
  correspondence between the ``top'' state for $K^+$ (whose local
  Alexander contribution is $1/2$, according to
  Figure~\ref{fig:AlexanderGradings}), with the same state, thought of
  as a state for the resolution.  Heegaard diagrams for a positive
  crossing and its destabilization. The black generator in the bottom
  row has local Alexander contribution $-1/2$, as can be seen by some
  handleslides.  The reader can easily find a Whitney disk connecting
  the black to the white state, which crosses $X$ with multiplicity
  $-1$, and which then corresponds to the given Kauffman state for the
  resolution (where we use the basepoints $B$ rather the ones marked
  by $X$'s).}
\end{figure}

Having established the behaviour of the Alexander grading, we find it
convenient to work with the algebraic, rather than the Maslov,
grading, since that gives us the added flexibility to
specialize to $U_i=1$
and then isotope across
the $O_i$.

\end{proof}

\section{Proof of Theorem~\ref{thm:Calculate}.}
\label{sec:Proof}

We assemble now the ingredients to verify Theorem~\ref{thm:Calculate}.

We will need one more fact about the homomorphisms $\Zip_p$ and
$\Unzip_p$ defined in Section~\ref{sec:ExactSeq}.

Let ${\mathcal X}$ be a singular knot and $\Smooth$ the singular knot
obtained by smoothing one of the singular points $p$ of
$\Singularize$, and consider their planar diagrams.  Let
$\x_0(\Singularize)$ be the cycle representing the generator of the
$\RingM$-module $H_*(C'(\Singularize))$.

\begin{figure}[ht]
  \mbox{\vbox{\epsfbox{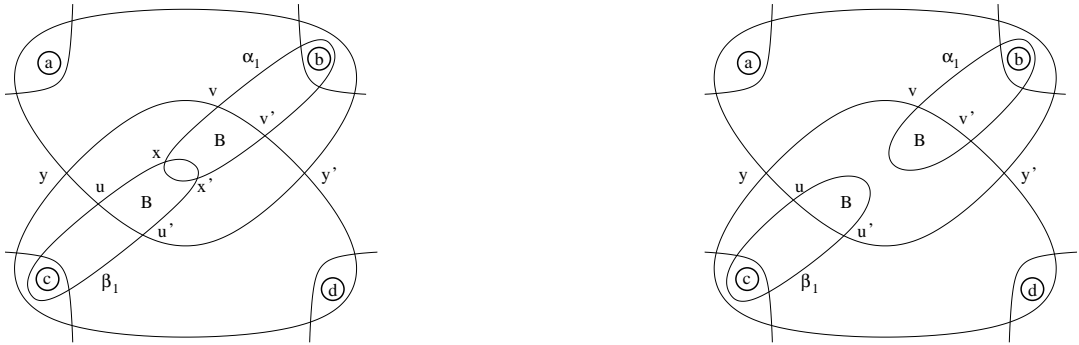}}}
\caption{\label{fig:ResolvedTwo}
  {\bf{Heegaard diagrams for the smoothing.}}  We have here two
  alternative Heegaard diagrams representing the smoothing at the
  distinguished crossing $p$.}
\end{figure}

To calculate $\Zip_p$, we use the Heegaard diagram for the smoothing
${\mathcal R}$ indicated in Figure~\ref{fig:ExactSequence} (using
basepoints $B$ as markings representing the smoothing $\Smooth$), and
also the left-hand diagram in Figure~\ref{fig:ResolvedTwo}.

Let $\y_0(\Smooth)$ denote the intersection point which, away from the
distinguished crossing agrees with the intersection point
$\x_0(\Smooth)$ of Proposition~\ref{prop:MinimalElement}, but which,
at the distinguished (smoothed) crossing consists of the pair of intersection 
points $x'$ and $y$ of Figure~\ref{fig:ResolvedTwo}.

We have the following analogue of Proposition~\ref{prop:MinimalElement}:

\begin{lemma}
  \label{lemma:MinimalElementSmoothing}
  The intersection point $\y_0({\mathcal R})$ is a cycle in the
  Heegaard diagram, and its image generates  $H_*(C'(\Smooth))$.
\end{lemma}

\begin{proof}
  Let $\y_1(\Smooth)$ be the intersection point which agrees with
  $\y_0(\Smooth)$ away from the distinguished crossing, where it exchanges
  $x'$ and $y$ for $u$ and $v$.

  Isotoping $\alpha_1$ and $\beta_1$ so that they are disjoint, as in
  the right-hand picture in Figure~\ref{fig:ResolvedTwo}, we obtain a
  new Heegaard diagram for the smoothing. For this new diagram, the
  proof of Proposition~\ref{prop:MinimalElement} applies: again, we
  can associate a coherent multi-cycle to each generator, and indeed,
  there is a generator with minimal algebraic grading among all
  generators associated to the empty multi-cycle, and it has has
  locally the form $u$ and $v$. In fact, this argument applies shows
  that of all the generators for original Heegaard diagram for the
  smoothing which do not contain $x$ or $x'$, the generator with
  minimal algebraic grading is $\y_1(\Smooth)$.
  
  Similarly, the proof of Proposition~\ref{prop:MinimalElement},
  applied now to the singularized diagram, shows that of the
  intersection points containing $x'$, there is a unique one with
  minimal algebraic grading, and in fact, it is $\y_0(\Smooth)$.

  But a small rectangle in the diagram, juxtaposed with a bigon with
  multiplicity $-1$, demonstrates that this generator has algebraic
  grading equal to that of $\y_1({\mathcal R})$. 
  Also the bigon
  from $x$ to $x'$ demonstrates that any generator contaning $x'$ has
  algebraic grading at least one greater than $\y_0(\Smooth)$.
  
  We have established that there are no generators with algebraic
  grading less than $\y_0(\Smooth)$, and hence it follows that
  $\y_0(\Smooth)$ is a cycle.

  Under the map induced by the isotopy, it is clear that
  $\y_1(\Smooth)$ is taken to a generator with minimal alegbraic
  grading in the new diagram, and hence it generates the homology
  $H_*(C'(\Smooth))$. Since it is homologous to $\y_0(\Smooth)$, it
  follows $\y_0(\Smooth)$ also generates this homology.
\end{proof}

\begin{prop}
  \label{prop:CalcMaps} Under the identifications
  $H_*(C'(\Singularize))\cong \AlgM(\Singularize)$ and
  $H_*(C'(\Smooth))\cong \AlgM({\mathcal R})$ (the first of which
  comes from Theorem~\ref{thm:CalcSingLink}, and the second is
  achieved with the help of
  Lemma~\ref{lemma:MinimalElementSmoothing}), the induced maps on
  homology $H_*(C'(\Singularize))\longrightarrow H_*(C'({\mathcal
  R}))$ and $H_*(C'(\Smooth))\longrightarrow H_*(C'({\mathcal X}))$
  induced by the unzip and zip maps of equation
  Equation~\eqref{eq:DefUnzipDiag} and~\eqref{eq:DefZipDiag}
  correspond to multiplication by $\pm 1$ and multiplication by $\pm
  (t\cm U_b-U_c)$ respectively (i.e. up to sign, they correspond to
  the maps $u_p$ and $z_p$ defined in Section~\ref{sec:Introduction}).
\end{prop}

\begin{figure}[ht]
  \mbox{\vbox{\epsfbox{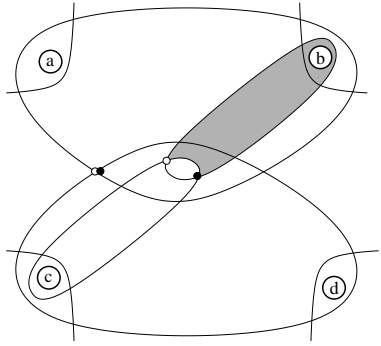}}}
\caption{\label{fig:ZipPicture}
  {\bf{Zip homomorphism.}} We have illustrated here one of the two
  bigons which represent the map
  $\Phi_B(\y_0(\Smooth))$.}
\end{figure}

\begin{proof}
  The map $\Unzip_p$ sends the subcomplex $X'$ of the mapping cone of
  $C'(\Singularize)$ (thought of as a mapping cone of $t\cm U_a+t\cm
  U_b-U_c-U_d\colon X' \longrightarrow X'$ isomorphically to the
  subcomplex $X'$ (thought of as a mapping cone of $\Phi_{A^+}\colon
  Y' \longrightarrow X'$. In particular, it carries the generator
  $\x_0(\Singularize)$ (in the notation of
  Proposition~\ref{prop:MinimalElement}) to the generator
  $\y_1(\Smooth)$ (in the notation of the proof of
  Lemma~\ref{lemma:MinimalElementSmoothing}). Thus, it follows at once
  that $\Unzip_p$ is induced by multiplication by $\pm 1$.
  
  The map $\Zip_p$, on the other hand, is induced from $\Phi_B$, i.e.
  it counts flowlines which cross one of the two regions marked by
  $B$.  We claim that there are exactly two Maslov index one, positive
  flowlines leaving the generator $\y_0(\Smooth)$ from
  Lemma~\ref{lemma:MinimalElementSmoothing}, crossing $B$, and landing
  at the canonical generator $\x_0(\Singularize)$. Both are
  bigons, one of which contributes multiplication by $\pm t\cm
  U_b$, and the other is multiplication by $\mp U_d$, cf. 
  Figure~\ref{fig:ZipPicture}.  The fact that these two flowlines
  contribute opposite signs follows from the fact that their
  difference is a difference of two boundary degenerations, one with
  $\alpha$-boundary, the other with $\beta$-boundary.
\end{proof}

\begin{proof}[Of Theorem~\ref{thm:Calculate}]
  Starting from a braid, we construct a planar diagram for it.  This
  gives us a cube of resolutions, as in
  Theorem~\ref{thm:SpectralSequence}.  The $E_1$ term of the
  associated Leray spectral sequence given by $$
  H_*\left(\CFKm(X_I(\Knot))\otimes \Big(\bigotimes_{s\in
  \Sing(X_I(\Knot))} \RingM/(t\cm U_a^{(s)}+t\cm
  U_b^{(s)}-U_c^{(s)}-U_d^{(s)})\Big)\right), $$ and its $d_1$
  differential is induced from the maps $\Zip_p$ and $\Unzip_p$ as
  defined above. These are calculated, up sign, in
  Proposition~\ref{prop:CalcMaps} above. The spectral sequence
  collapses after the $E_1$ stage, for the following reason. The
  differentials drop Maslov grading by one, and it is easy to see that
  they preserve the renormalized Alexander grading $A'$ as defined in
  Equation~\eqref{eq:RenormalizedAlexander}. It follows that the
  differential drops the algebraic grading by one, as well. But
  according to Theorem~\ref{thm:CalcSingLink} the algebraic grading at
  each stage is proportional to the filtration level, and hence the
  spectral sequence collapses.
  
  It follows that
  $$H_*(\CnewM(\Knot)\otimes_{\Z[t]}\Z[[t]])
  \cong H_*(\CFKmT(\Knot)\otimes_{\Z[t]}\Z[[t]]).$$
  A similar reasoning gives
  $$H_*(\CnewM(\Knot)\otimes_{\RingM} \RingM[[t]]/U_0)
  \cong H_*(\CFKaT(\Knot)\otimes_{\Z[t]}\Z[[t]]).$$

  The identifications stated in Theorem~\ref{thm:Calculate} follow
  from these isomorphisms, together with
  Lemma~\ref{lemma:TwistedCoefficients}, once the signs have been taken
  care of.
  
  The signs in the introduction are set up so that $D^2=0$ in the cube
  of resolutions (i.e. on the chain level, rather than merely on its
  $E_2$ term). In fact, it is easy to see that they are uniquely
  determined, up to an overall isomorphism
  of the chain complex, by this property.
\end{proof}

\section{Grid diagrams and singular links.}
\label{sec:Grids}

The skein exact sequence involving singular links can be readily seen
from the point of view of grid diagrams, following~\cite{MOS}, \cite{MOST}.

We set up some of the terminology first.

A {\em planar grid diagram}~$G$ lies in a square
grid on the the plane with $n \times n$ squares, each
of which is decorated by $X$, $O$, or nothing, with the following
rules:
\begin{itemize}
\item {every square is decorated by exactly one $X$, exactly one $O$,
    or nothing}
\item {every row contains exactly one $X$ and one $O$;}
\item {every column contains exactly one $X$ and one $O$.}
\end {itemize}
A planar grid diagram can be
transferred to the torus, by identifying the opposite sides of the
$n\times n$ grid above. In this manner, we obtain a Heegaard diagram
for an oriented link in the three-sphere, in the sense of
Section~\ref{sec:KnotFloerHomology}.

Given a planar grid diagram, we can in fact construct braid as
follows, cf.~\cite{BirmanMenasco}, \cite{Cromwell}. Draw oriented
segments in each row starting at the corresponding $X$ and ending at
the corresponding $O$, and then draw oriented from segments in each
column starting at the corresponding $O$ and ending at the
corresponding $X$, provided that $X$ is above the $O$. Otherwise, we
draw an ``outgoing'' arc starting at $O$, and an ``incoming'' one
ending at $X$. At each crossing of the horizontal and vertical arcs,
we take the convention that the vertical arc is an overcrossing.
Indeed, this braid is naturally associated to the picture of the grid
diagram, as drawn on the cylinder, where the top and the bottom sides
of the planar grid diagram are identified.

Certain squares in the grid diagram correspond to crossings in the
projection; these squares are called {\em crossing squares}.  An $X$
in the same row or column as a crossing square is called a {\em
  crossing $X$}; all others are called {\em non-crossing $X$'s}.
Finally, each crossing square has a distinguished corner which is also
contained in two of the crossing $X$'s.

For us, braids will always be {\em decorated}: this means, we place a marker
on one of the leftmost edges of our braid.

\begin{defn}
\label{def:Special}
A planar grid diagram for a decorated braid is called {\em special} if
the following conditions are satisfied:
\begin{enumerate}
\item each vertical arc crosses over at
most one horizontal arc (or, equivalently, each row and column
contains at most one crossing square)
\item no two crossing squares share a corner
\item each crossing square shares two sides with squares marked by
  $X$'s; in this case we call the corner shared by the two $X$'s a
  {\em crossing corner}
\item 
  \label{item:KillRectangles}
  any rectangle which has the property that two of its corners are
  crossing corners also has an $X$ in its interior.
\item 
  the segment above the decoration corresponds to the leftmost
  vertical arc in the grid diagram.
\end{enumerate}
\end{defn}

It is easy to construct special grid diagrams for a decorated braid,
after sufficiently many stabilizations.

\begin{figure}[ht]
\mbox{\vbox{\epsfbox{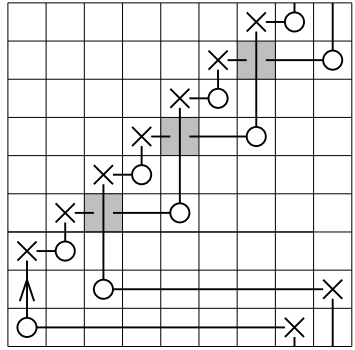}}}
\caption{\label{fig:SpecialBraid}
  {\bf{Special braid.}}  We illustrate here a special grid diagram for
  a braid form of the trefoil knot. The crossing squares are shaded.}
\end{figure}

A braid diagram for a singular link is a diagram where
\begin{itemize}
\item every square is decorated by exactly one or two $X$'s, 
    exactly one or two $O$'s, or nothing
\item the total number of $X$'s in a row or column equals the total
  number of $O$'s, which in turn equals one or two.
\end {itemize}

Given such a diagram, it is straightforward to construct the
corresponding singular link as before, with the understanding that the
squares marked by two $X$'s correspond to the singular points, cf.
Figure~\ref{fig:SingularBraid} below.

\begin{figure}[ht]
\mbox{\vbox{\epsfbox{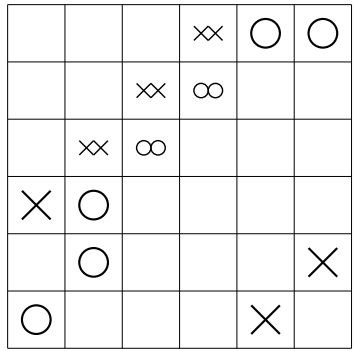}}}
\caption{\label{fig:SingularBraid}
  {\bf{Singular braid.}}  A grid diagram for the singularization of
  the trefoil from Figure~\ref{fig:SpecialBraid}.}
\end{figure}

As in~\cite{MOS}, the chain complexes for $\CFKm$ and $\CFKa$ count
only rectangles, cf. also~\cite{SarkarWang}; see~\cite{MOST} for the
sign refinement.

We now turn to a grid diagram proof of Theorem~\ref{thm:SkeinExactSequence}.

Suppose we have four squares in a grid diagram which meet at a corner
$c$, and whose complement contains all the $O_i$ and $n-2$ of the
$X$'s, marked $\Xs_0$.  Next, mark the upper left and lower right
squares by $A$ and mark the upper right and lower left square by $B$
and (we call these pairs of squares $\As$ resp. $\Bs$). We can form
alternative grid diagrams $G_A$ resp. $G_B$, both of which use the
same set $\Os$, and using $\Xs=\As\cup \Xs_0$ resp. $\Bs\cup \Xs_0$.
Let $O_a$ and $O_b$ be the two $O$'s in the columns through $A$ and
$B$, and let $O_c$ and $O_d$ be the two $O$'s in the rows through $A$
and $B$. 

\begin{figure}[ht]
\mbox{\vbox{\epsfbox{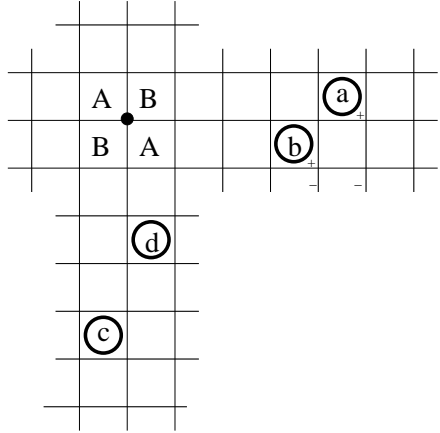}}}
\caption{\label{fig:AB}
{\bf{Squares marked $A$ and $B$.}}
We have also included positive and negative markers, denoted here by $+$ and $-$.
Note that the squares marked $O_c$ and $O_d$ (and indeed another square in their respective column)
might or might not contain additional markings, depending on whether or not $O_c$ and $O_d$ belong
to crossings in the diagram.}
\end{figure}

Clearly, $\CFKm(G_B)$ has a subcomplex $X$ consisting of configurations
which contain the corner $p$, and a quotient complex $Y$. Moreover,
$\CFKm(G_A)$ has $Y$ as a subcomplex, with quotient $X$. Thus, 
$\CFKm(G_B)$ can be thought of as the mapping cone of the map
$$\Phi_A\colon Y \longrightarrow X$$
gotten by counting rectangles which contain exactly one of the squares
marked by $B$, 
and $\CFKm(G_A)$ is
the mapping cone of the map
$$\Phi_B\colon X \longrightarrow Y,$$
defined by counting rectangles
which contain exactly one of the squares marked by $A$.
More precisely,
$$\Phi_B(\x)=\sum_{\y\in{\mathfrak S}} \sum_{\{r\in\EmptyRect(\x,\y)
  \big| X_i(r)=0~~~~~~~~~~\forall i\}} \sign(r)\cm 
U_1^{O_1(r)}\cm...\cm U_n^{O_n(r)}\cm \y, $$
where $\sign$ denotes the sign assignment $\pm 1$ to each rectangle.

We have the following analogue of Lemma~\ref{lemma:CalcComposite}:

\begin{lemma}
  \label{lemma:CompositeMap}
The composite map $\Phi_A\circ \Phi_B$ is equal to multiplication by
$$U_a+U_b-U_c-U_d.$$
\end{lemma}

\begin{proof}
  Given some generator $\x$ which contains the corner $c$, The
  double-composite $\Phi_A\circ\Phi_B$ counts (with suitable signs and
  $U_i$-powers) the four annuli (two horizontal and two vertical)
  which contain $c$ on their boundary, cf.
  Figure~\ref{fig:CalcComposite}.  The horizontal annuli contribute
  $U_a$ and $U_b$, while the vertical ones contribute $-U_c$ and
  $-U_d$. To verify the signs, recall~\cite{MOST} that $\sign$ is
  determined so that if $r_1$ and $r_2$ are two empty rectangles whose
  composite is an annulus, then $\sign(r_1)\sign(r_2)=+1$ if the
  annulus is horizontal, and $-1$ if it is vertical.
\end{proof}

\begin{figure}[ht]
\mbox{\vbox{\epsfbox{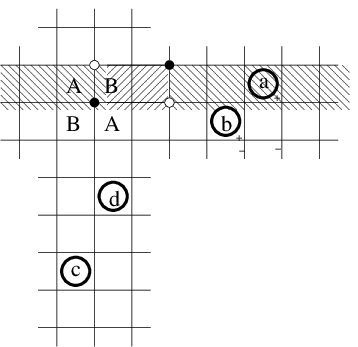}}}
\caption{\label{fig:CalcCompositeGrid}
{\bf{Composite $\Phi_A\circ\Phi_B$.}}}
\end{figure}

As before, we consider the following two bigraded complexes $P$ and $N$:
\begin{eqnarray*}
\begin{CD}
  X @>{\Phi_B}>> Y \\
  @V{U_a+U_b-U_c-U_d}VV @VV{\Phi_A}V \\
  X  @>{=}>> X
\end{CD}
&\text{{\hskip.5cm}and{\hskip.5cm}}&
\begin{CD}
  X @>{=}>> X \\
  @V{\Phi_B}VV @VV{U_a+U_b-U_c-U_d}V \\
  Y  @>{\Phi_A}>> X.
\end{CD}
\end{eqnarray*}
Note that, by contrast with the situation in Section~\ref{sec:ExactSeq},
there is no need for homotopies $\Phi_{AB}$ and $\Phi_{BA}$: they both vanish
(since our diagram is special in the sense of Definition~\ref{def:Special}),
and consequently, the above constructions do indeed give complexes.

The first of these is clearly
quasi-isomorphic to $\CFKm(G_A)$, while the second is quasi-isomorphic
to $\CFKm(G_B)$. This can be seen by considering the filtration
induced by the vertical coordinate, where the complex is expressed as
the extension of a subcomplex with trivial homology by a quotient
complex which is $\CFKm(G_A)$ (compare Section~\ref{sec:ExactSeq}). Alternatively, if we consider the the
filtration by the horizontal coordinate, where the first complex has a
subcomplex quai-isomorphic to $C_B$, with quotient complex the mapping
cone of
$$U_a+U_b-U_c-U_d \colon X \longrightarrow X;$$
or the second, which has a subcomplex quasi-isomorphic to the mapping cone of
$$U_a+U_b-U_c-U_d \colon X \longrightarrow X$$
and quotient complex quasi-isomorphic to $C_B$.

Using either observation, we obtain a long exact sequence
\begin{equation}
  \label{eq:SkeinExactSequence}
  \begin{CD}
    ...@>>>  H_*(\frac{X}{U_a+U_b- U_c - U_d}) @>>> H_*(C_B) @>>> H_*(C_A) @>>>
  \end{CD}
\end{equation}

There are now two cases: either $G_A$ and $G_B$ represent $\Ko$ and $\Km$
respectively, or $G_A$ and $G_B$ represent $\Kp$ and $\Ko$ respectively.

Indeed, one can iterate this procedure to give an alternate proof of 
Theorem~\ref{thm:SpectralSequence}, except when one uses special grid
diagrams, then all the higher homotopies present in the construction of
the chain complex vanish identically.

To give the construction with twisted coefficients, we must modify the construction of twisted
coefficients a little. Specifically, we associate it now not simply to a positive marking 
${\mathbb P}=\{P_1,...,P_m\}$ as before, but to a collection of positive and negative markings
${\mathbb P}=\{P_1,...,P_a\}$, ${\mathbb Q}=\{Q_1,...,Q_b\}$, with coefficients over $\Z[t,t^{-1}]$,
with the understanding that a disk $\phi$ contributes a factor of
$t^{P(\phi)-Q(\phi)}$.  This can be done only with coefficients in $\Z[t,t^{-1}]$ (as opposed to $\Z[t]$).

Now, for our special braid diagram, we must include markings as
follows.  Out of each crossing, there is a pair of $O$'s, $O_a$ and
$O_b$; if the crossing is positive resp. negative, then $O_b$ occupies
the row immediately above resp. below $O_a$.  Write $O_1=O_a$ resp.
$O_b$ if the crossing is positive resp. negative, and let $O_2$ be the
other one of $O_a$ or $O_b$ (so that $O_1$ is immediately below
$O_2$).  We mark $O_1$ with a positive marker, and the rectangle
immediately below it with a negative marker.  We mark $O_2$ with a
positive marker, and the rectangle immediately below it with nothing,
and the rectangle immediately below that with a negative marker.
The markings are indicated in Figure~\ref{fig:AB} (the $P$'s with $+$,
and $Q$'s with $-$).

With these markings, Lemma~\ref{lemma:CompositeMap} can be modified:
the double composite corresponds to multiplication by $t\cm U_a+t\cm
U_b-U_c-U_d$. This gives an alternate proof of
Theorem~\ref{thm:SpectralSequence}, with coefficients in
$\Z[t,t^{-1}]$, entirely within the realm of grid diagrams

The grid diagram on its face value is too unwieldly to verify the
results from Section~\ref{sec:CalcSingLink}.  However, appealing to
those results, one can also calculate the zip and unzip homomorphisms,
as in Proposition~\ref{prop:CalcMaps}, to complete an alternate proof
of Theorem~\ref{thm:Calculate}.  Specifically, we have the following:

\begin{defn}
  Let $\Lambda$ be the chain represented by the intersection points
  which occupies the lower left corner of each square marked with one
  or two $X$.
\end{defn}

\begin{lemma}
  The element $\Lambda$ is a cycle, which represents the generator of
  $H_*(C')$ over $\RingM$.
\end{lemma}

\begin{proof}
  $\Lambda$ is a cycle since each rectangle leaving it must contain
  one or two $X$ (compare~\cite{Transverse}). In fact, it is easy to
  see that $\Lambda$ has maximal Alexander grading among all
  generators. Thus, according to Theorem~\ref{thm:CalcSingLink}, it
  generates $H_*(C'(S))$.
\end{proof}

\begin{figure}[ht]
  \mbox{\vbox{\epsfbox{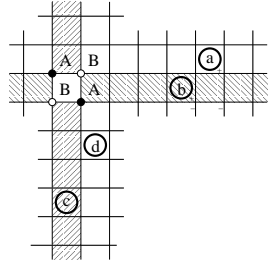}}}
\caption{\label{fig:CalcZip}
  {\bf{Zip homomorphism.}} We have illustrated here the two
        rectangles which contribute to $\Zip_p(\Lambda_r)$.}
\end{figure}

\begin{lemma}
  Let $\Lambda_r$ and $\Lambda_x$ be the canonical generators for
  the grid diagrams for ${\mathcal R}$ and ${\mathcal X}$ respectively.
  Then, $\Phi_A(\Lambda_r)=(t\cm U_a-U_d) \Lambda_x$.
\end{lemma}

\begin{proof}
  This follows from the fact that are precisely two Maslov index one
  squares leaving $\Lambda_r$ and ending in $\Lambda_x$, and they are
  the ones illustrated in Figure~\ref{fig:CalcZip}.
\end{proof}

The above lemma gives the needed geometric ingredient for the alternate
proof of Proposition~\ref{prop:CalcMaps}. 

\section{Examples}

\subsection{Illustrations of Theorem~\ref{thm:CalcSingLink}.}

\subsubsection{The singularized trefoil}
\label{subsubsec:Trefoil}
\begin{figure}[ht]
\mbox{\vbox{\epsfbox{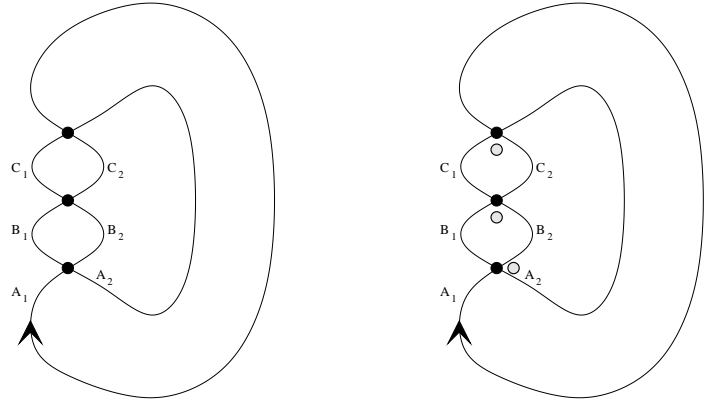}}}
\caption{\label{fig:SingTref}
  {\bf{Singularized trefoil.}}  We show the labelings for the
  singularized trefoil. The arrow indicates the orientation, and also
  the distinguished edge.  At the right we have illustrated the only
  possible type of Kauffman state (by gray dots). Of course, there are
  really four Kauffman states, according to which of the two types of 
  type $D$ Kauffman we use.}
\end{figure}

We begin with an illustration of the techniques from
Section~\ref{sec:CalcSingLink}, considering the the two-braid which is
the singularized trefoil pictured in
Figure~\ref{fig:SingTref}.

Label edges $A_1$, $A_2$, $B_1$, $B_2$, $C_1$, $C_2$, and the final half edge by $0$, as shown in that figure.
We have three linear relations coming from the Koszul complex:
\begin{eqnarray*}
  A_1+A_2&=&t B_1 + t B_2 \\
  B_1+B_2 &=& t C_1 + t C_2 \\
  C_1 + C_2&=& t A_2,
\end{eqnarray*}
one relation coming from the point set consisting of all vertices:
$$
A_1 = 0,
$$
and three quadratic relations
\begin{eqnarray*}
  A_1 A_2 &=& t^2 B_1 B_2 \\
  B_1 B_2 &=& t^2 C_1 C_2 \\
  C_1 C_2 &=& 0
\end{eqnarray*}
Substitute $A_1=0$ into the first linear equation; then substitute the second two to get 
$A_2=t^3 A_2$, from which it follows that $A_2=0$, as well.
We thus get the following equivalent set of equations:
\begin{eqnarray*}
  B_1+B_2&=&0 \\
  C_1+C_2&=&0 \\
  B_1 B_2 &=& 0 \\
  C_1 C_2 &=& 0.
\end{eqnarray*}
Eliminating $B_1$ and $C_2$, and writing $B$ and $C$ for $B_2$ and $C_2$
respectively, we see that the ring is given by
$$\AlgM(S)/U_0\cong \Z[B,C]/\{B^2=0, C^2=0\}.$$  This ring is four-dimensional over $\Z[t,t^{-1}]$,
generated by  $1$, $B$, $C$, and
$BC$. This is consistent with the state sum formula, according to
which there is only one Kauffman state, and it is the one illustrated
on the left of Figure~\ref{fig:SingTref}.

\subsubsection{The singularization of $4_1$} Consider the singularization of
the four-crossing number. The corresponding algebra in this case is somewhat
more interesting. 

\begin{figure}[ht]
\mbox{\vbox{\epsfbox{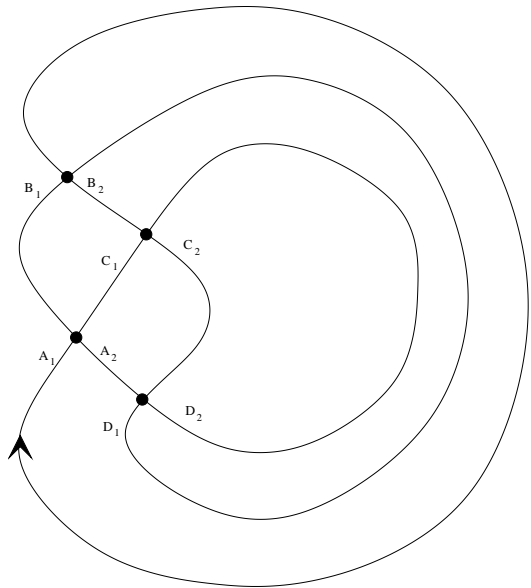}}}
\caption{\label{fig:SingFourOne}
  {\bf{Singularized $4_1$.}}  We show the labelings for the
  singularization of $4_1$.}
\end{figure}

Labeling the variables corresponding to the edges $\{A_1,A_2,B_1,B_2,C_1,C_2,D_1,D_2\}$ as indicated in Figure~\ref{fig:SingFourOne},
we obtain four linear relations:
\begin{eqnarray*}
  A_2&=&tB_1+tC_1 \\
  C_1+C_2&=& tB_2+tD_2 \\
  D_1+D_2&=&tA_2+tC_2 \\
  B_1+B_2&=&tD_1.
\end{eqnarray*}
The point set consisting of all vertices gives the relation 
$$A_1=0.$$
There are four quadratic relations coming from the four vertices
\begin{eqnarray*}
A_1 A_2 &=& t^2 B_1 C_1 \\
B_1 B_2&=&0 \\
C_1 C_2 &=& t^2 B_2 D_2 \\
D_1 D_2 &=& t^2 A_2 C_2,
\end{eqnarray*}
and one more quadratic relation which will come into play, coming from the small innermost cycle
$$C_1 D_1 = t^4 A_2 B_2.$$

The linear equations can be used to eliminate $A_2$, $B_1$, $C_2$, and $D_2$ writing them in terms of $B_1$, $C_1$, and $D_1$:
\begin{eqnarray*}
A_2&=&D_1 t^2-B_2 t+C_1 t \\
B_1&=&D_1
   t-B_2 \\
C_2&=&-C_1+B_2 t-\frac{t   \left(D_1 t^2+C_1 t+D_1
   t+D_1\right)}{t+1} \\
D_2&=&-\frac{D_1
   t^2+C_1 t+D_1
   t+D_1}{t+1} 
\end{eqnarray*}
Substituting these back into the quadratic equations, we see that the quadratic part of the algebra is one-dimensional, generated
by $B_2 D_1$; indeed,
\begin{eqnarray*}
B_2^2 &=& B_2 D_1 t \\
B_2 C_1 &=& 0 \\
C_1^2 &=& B_2 D_1 t^2 \\
C_1 D_1 &=& 0 \\
D_1^2 &=&-\frac{{B_2 D_1} t^4}{t^4+t^3+t^2+t+1} 
\end{eqnarray*}
It is also easy to see that the degree three part of the algebra is trivial. Thus, we have verified that the is five dimensional over $\Z[t^{-1},t]]$,
generated by $1$, $B_2$, $C_1$, $D_2$, and $B_2 D_1$, consistent with a straightforward count of Kauffman states.

\subsection{Calculations using the cube of resolutions}

We illustrate Theorem~\ref{thm:SkeinExactSequence}, as a calculational device. We consider three-crossing examples. 
Consider the three-crossing projection of the right-handed trefoil.

\begin{figure}[ht]
\mbox{\vbox{\epsfbox{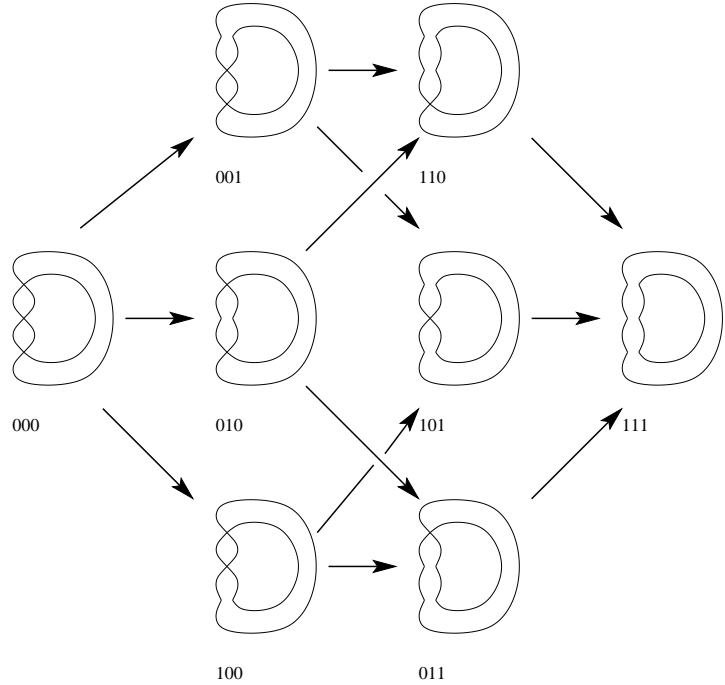}}}
\caption{\label{fig:CubeResTref}
  {\bf{Cube of resolutions for the right-handed trefoil.}}  We have
  indicated the resolution type with a string of zeros and ones.
  Differentials in the cube of resolutions are indicated by arrows.
  The horiztontal cooordinate is the ``algebraic grading''.}
\end{figure}

The cube of resolutions has a natural horizontal grading (each
resolution is indexed by three elements in $\{0,1\}$; the horizontal
grading is the sum of the three numbers). The Alexander polynomial of
the singular link in the leftmost column is $T^{-1}+2+T$; the
Alexander polynomial of each singular link in the next column is
$-(T^{-1/2}+T^{1/2})$, but these are shifted by $T^{1/2}$ in the cube of
resolutions; the Alexander polynomial of each singular link in the
next column is $1$, and this is shifted by $T$; and the Alexander
polynomial in the final column is $0$.  In effect, we have decomposed
the Alexander polynomial of the trefoil as:
$$T^{-1}-1+T = (T^{-1}+2+T) +3 (-1-T) + 3 (T) + 0.$$

The cube of resolutions has a corresponding splitting by Alexander
gradings. There is a single generator in Alexander grading $-1$, and
it is supported in the $000$ singularization, where it has Maslov
grading equal to $-2$.

Next, we consider the part in Alexander grading $0$. There are two
generators coming from the $000$ singularization which, in the
notation from Subsection~\ref{subsubsec:Trefoil}, is generated by the
monomials $B$ and $C$; these both have Maslov grading equal to zero.
There are three further generators in Alexander grading $0$, each of
which has Maslov grading $-1$, coming from degree one monomials
belonging to the FLoer homologies of the singularizations $001$,
$010$, and $100$. There are no further generators with Alexander
grading zero, so it suffices to show that the differential carries the
subspace generated by $B$ and $C$ injectively into this
three-dimensional space. It is easy to see that $B$ is annihilated by
its map into the $100$ resolution, although $C$ injects there; while
$C$ is annihilated by its map into $001$, while $B$ injects there.  It
follows easily that the summand of knot Floer homology in Alexander
grading equal to zero is one-dimensional, supported in Maslov grading
$-1$.

A similar calculation verifies that the remaining Floer homology group
is supported in Alexander grading $+1$ and Maslov grading $0$,
compare~\cite{Knots}.

\bibliographystyle{plain}
\bibliography{biblio}

\end{document}